\documentclass[11pt]{article}
\usepackage[a4paper,margin=2.5cm]{geometry}
\usepackage[comma, round]{natbib}
\usepackage{float}
\usepackage{amssymb}
\usepackage{amsthm}
\usepackage{amsmath}
\usepackage{titletoc}
\usepackage{mathrsfs}
\usepackage[colorlinks, allcolors=blue]{hyperref}
\usepackage{amsfonts}
\usepackage{graphicx}
\usepackage[linesnumbered, ruled, vlined]{algorithm2e}
\usepackage{esint}
\usepackage{bbm}
\usepackage{mathtools}
\usepackage[shortlabels]{enumitem}
\usepackage{authblk}
\usepackage[capitalize,noabbrev]{cleveref}

\numberwithin{equation}{section}
\numberwithin{figure}{section}

\theoremstyle{plain}
\newtheorem{thm}{Theorem}[section]
\newtheorem{cor}[thm]{Corollary}
\newtheorem{lem}[thm]{Lemma}
\crefname{lem}{Lemma}{Lemmas}
\newtheorem{prop}[thm]{Proposition}

\theoremstyle{definition}
\newtheorem{defn}{Definition}[section]
\newtheorem{asmp}{Assumption}[section]
\newtheorem{rmk}[defn]{Remark}

\newcommand{\1}{\mathbbm{1}}
\newcommand{\AW}{\mathcal{AW}}
\newcommand{\W}{\mathcal{W}}

\newcommand{\D}{\mathcal{D}}

\newcommand{\md}{\mathop{}\mathopen\mathrm{d}}
\newcommand{\E}{E}

\renewcommand{\P}{P}
\renewcommand{\L}{\mathcal{L}}

\newcommand{\Id}{\operatorname{Id}}

\renewcommand{\c}{\mathfrak{c}}
\newcommand{\bc}{\mathfrak{bc}}

\newcommand{\bt}{\boldsymbol{\tau}}
\newcommand{\bd}{\boldsymbol{\delta}}

\title{Stability, discretization, and optimality of continuous-time causal optimal transport}
\author[1]{Yifan Jiang}
\author[2]{Fang Rui Lim}
\affil[1]{Department of Mathematics, Imperial College London}
\affil[2]{Department of Mathematics, University of Michigan, Ann Arbor}
\date{}

\begin{document}
\maketitle

\begin{abstract}
	We prove stability of continuous-time causal and bicausal optimal transport problems when both marginals converge in adapted Wasserstein distance along common vanishing time shifts.
	The marginal laws, filtrations, and costs may all vary, while the limit remains the original causal or bicausal transport problem.
	This allows for time discretization and the approximation of transport problems between  non-degenerate scalar fractional stochastic differential equations  by their Wong--Zakai schemes.
	We further establish sufficient conditions so that the synchronous coupling is optimal for the bicausal transport problem between scalar fractional stochastic differential equations under continuous, nonnegative submodular costs of polynomial growth.
	This result covers all Hurst parameters \(H\in(0,1)\) up to an additional  order-preserving commutator condition on the  corresponding coefficients for the rough regime \(H\in(0,1/2)\).
\end{abstract}

\begin{description}
	\item[Mathematics Subject Classification (2020):]  	49Q22, 60G05, 60G22
	\item[Keywords:]  causal optimal transport, adapted Wasserstein distance, stability, fractional SDEs, Wong--Zakai approximation.
\end{description}

\section{Introduction}

A central problem in stochastic analysis is to understand when a continuous-time object can be interpreted as a limit of discrete-time approximations.
This theme underlies invariance principles, martingale approximations, rough path theory, and numerical schemes for stochastic equations.
Such limiting procedures often require a suitable rescaling or a controlled relaxation of the discrete objects.
The classical example is Donsker's invariance principle \citep{donsker51Invariance}.
We consider a standard Brownian motion \(B\) and its polygonal approximation \(B^{n}\) given by
\begin{equation}
	\label{eqn-donsker}
	B^n_t := B_{\lfloor nt\rfloor/n}+(nt -\lfloor nt\rfloor) (B_{(\lfloor nt\rfloor+1)/n}-B_{\lfloor nt\rfloor/n})\qquad t\in[0,T].
\end{equation}
Then \(B^n\) converges in distribution to  \(B\) in the uniform topology on \(C([0,T];\mathbb{R})\); in fact, the convergence holds almost surely in this special case.
Thus, the interpolated process \(B^n\) approximates the path law of \(B\).
However, weak convergence does not record how information is revealed: the interpolated Brownian motion releases information only at grid times, whereas the Brownian filtration evolves continuously.
This distinction is fundamental for optimization under adaptedness constraints, but remains invisible to weak convergence.

The need to retain this information flow led to the notions of adapted distributions and the adapted weak topology \citep{hoover84Adapteda,aldous81Weak,hellwig96Sequential}.
More recently, these notions have been linked to \emph{causal} optimal transport, an optimal transport problem on path space with an adaptedness constraint on the coupling; see \citet{lassalle18Causal,acciaio20Causal,backhoff-veraguas20All,bartl24Wasserstein}.
Let \(X,Y\) be stochastic processes equipped with filtrations \(\mathbf{F},\mathbf{G}\), which need not be their natural filtrations.
A coupling between \((X,\mathbf{F})\) and \((Y,\mathbf{G})\) is \emph{causal} if the target up to time \(t\) cannot use future information from the source, and \emph{bicausal} if this holds in both directions.
Formally, the (bi)causal optimal transport problem is defined as
\begin{equation}
	\label{eqn-intro-transport}
	\mathcal{T}_{\bullet}((X,\mathbf{F}),(Y,\mathbf{G});c) :=\inf_{\pi\in\Pi_{\bullet}(\mathbf{F},\mathbf{G})}\E_{\pi}[c(X,Y)], \qquad \bullet\in\{\c,\bc\},
\end{equation}
where \(\Pi_{\c}(\mathbf{F},\mathbf{G})\) and \(\Pi_{\bc}(\mathbf{F},\mathbf{G})\) denote the corresponding sets of causal and bicausal couplings.
Precise definitions are deferred to \cref{sec-prelim}.
The so-called adapted Wasserstein distance is given by
\begin{equation*}
	\AW((X,\mathbf{F}),(Y,\mathbf{G}))=\mathcal{T}_{\bc}((X,\mathbf{F}),(Y,\mathbf{G}); d),
\end{equation*}
where \(d\) is a metric on the path space.
In the seminal paper \citet{backhoff-veraguas20All}, it is shown that the adapted Wasserstein distance metrizes the adapted weak topology in discrete time.

Much of the structural, stability, and computational theory of causal optimal transport has been developed in discrete time, where the information flow is indexed by finitely many stages \citep{backhoff-veraguas17Causal,backhoff-veraguas20All,eckstein24Computational}.
Continuous-time causal optimal transport, on the other hand, is often more rigid and involves more delicate analysis.
In particular, if \(B^n\) from \eqref{eqn-donsker} and its Brownian limit \(B\) are equipped with their natural filtrations, \citet{hoover84Adapteda} established that \(B^{n}\) converges to \(B\) in the adapted weak topology.
However, the adapted Wasserstein distance between them does not converge to zero because their information is revealed at different times and the independent coupling is the only bicausal coupling.
This is in sharp contrast to the discrete-time result and breaks the bridge between continuous- and discrete-time frameworks.

To address this issue, we introduce convergence in adapted Wasserstein distance along a sequence of vanishing time shifts \(\bt=(\tau^{n})_{n\geq 1}\).
For an increasing map \(\tau:[0,T]\to[0,T]\) with \(\tau(t)\geq t\), we write \(\mathbf{F}\circ\tau=(\mathcal{F}_{\tau(t)})_{t\in[0,T]}\).
This time shift perturbs information flows, not paths.
If \(\|\tau^n-\Id\|_{\infty}\to 0\), we say that \((X^n,\mathbf{F}^n)\) converges to \((X,\mathbf{F})\) in adapted Wasserstein distance along \(\bt=(\tau^n)_{n\geq 1}\) when
\begin{equation}
	\label{eqn-intro-shift-conv}
	\AW((X^n,\mathbf{F}^n\circ\tau^n),(X,\mathbf{F}\circ\tau^n))\longrightarrow0.
\end{equation}

This brings us back to Donsker's theorem.
We take \(X^{n}=B^{n}\), \(X=B\), and \(\tau^{n}:t\mapsto \frac{T}{n} (\lceil \frac{t(n-1)}{T} \rceil+1)\), which shifts all times to the right end of the grid interval.
The synchronous coupling of \(B^n\) and \(B\) is bicausal after the common shift \(\tau^{n}\), and uniform continuity of Brownian paths gives the convergence in \eqref{eqn-intro-shift-conv}.
Thus, convergence along vanishing time shifts accommodates this basic polygonal approximation even though strict adapted Wasserstein convergence fails.

Our first main  contribution shows that this shifted convergence is sufficient for stability of \eqref{eqn-intro-transport} in continuous time.
If both marginals converge along the same sequence of vanishing time shifts, the limiting filtrations are right-continuous after completion, and the costs satisfy standard convergence and uniform-integrability conditions, \cref{thm-stab} yields
\begin{equation}
	\label{eqn-intro-stability}
	\mathcal{T}_{\bullet}\bigl((X^n,\mathbf{F}^n\circ\tau^n),(Y^n,\mathbf{G}^n\circ\tau^n);c_n\bigr) \longrightarrow \mathcal{T}_{\bullet}((X,\mathbf{F}),(Y,\mathbf{G});c), \qquad \bullet\in\{\c,\bc\}.
\end{equation}
Thus, the marginal laws, filtrations, and costs may all vary, while the limit on the right remains the original causal or bicausal problem.
The proof combines convergence of the shifted causal coupling sets with a quantitative gluing argument that exchanges couplings between the original filtrations and the shifted ones.

Our second contribution studies the (bi)causal optimal transport between scalar fractional SDEs given by
\begin{equation}
	\label{eqn-intro-fsde}
	\left\{
	\begin{aligned}
		\md X_{t} & = b(t,X_{t}) \md t + \sigma(t,X_{t}) \circ \md Z_{t}\text{ with } X_{0}=x_{0},                         \\
		\md Y_{t} & = \bar{b}(t,Y_{t}) \md t + \bar{\sigma}(t,Y_{t}) \circ \md \bar{Z}_{t} \;\; \text{ with } Y_{0}=y_{0},
	\end{aligned}
	\right.
\end{equation}
where \(Z\) and \(\bar{Z}\) can be a standard Brownian motion or a fractional Brownian motion.
The solution here is understood in the Doss--Sussmann sense \citep{doss77Liens,sussmann78Gap}, and  all processes are equipped with their natural filtrations.
In Section \ref{sec-disc}, we construct a discrete approximation \(X^{n}\) of \(X\) by replacing \(Z\) in \eqref{eqn-intro-fsde} with the piecewise linear approximation \(Z^{n}\) on a dyadic grid \(D_{n}\).
We show that \(X^{n}\) converges to \(X\) in adapted Wasserstein distance along a sequence of vanishing time shifts \(\bt\).
Moreover, applying the stability result yields Theorem \ref{thm-approx}, which expresses the original continuous-time causal transport problem as the limit of a sequence of discrete-time problems:
\begin{equation}
	\label{eqn-intro-approx}
	\mathcal{T}_{\bullet}(X,Y;c)=\lim_{n\to\infty}\mathcal{T}_{\bullet}(\hat{X}^n,\hat{Y}^n;\hat c_n), \qquad \bullet\in\{\c,\bc\},
\end{equation}
where \(\hat{X}^{n}=(X^{n}_{t})_{t\in D_{n}}\) is the discrete-time process extracted from \(X^{n}\) on the grid \(D_{n}\), \(\hat{Y}^{n}\) is defined analogously, and \(\hat c_n\) denotes the cost of the interpolated grid paths.
From a numerical perspective, one can further apply  existing discrete-time computational methods in \citet{backhoff-veraguas17Causal,eckstein24Computational}.

Our third contribution further exploits the above discrete approximation to identify an optimal bicausal coupling for scalar fractional SDEs.
For paths \(x\) and \(y\), we write \(x\preceq y\) if \(x(t)\leq y(t)\) for any \(t\in[0,T]\).
We assume that \(c\) is a continuous submodular cost on the path space in the sense that
\begin{equation*}
	c(x,y)+c(x',y')\leq c(x,y')+c(x',y) \text{ for any paths } x\preceq x' \text{ and }y\preceq y'.
\end{equation*}
For \((X,Y)\) as in \eqref{eqn-intro-fsde}, the \emph{synchronous coupling} uses the same underlying Brownian motion \(B\) in both Molchan--Golosov representations \eqref{eqn-Molchan-Golosov} of \(Z\) and \(\bar{Z}\).
We show in Theorem \ref{thm-opt} that this  coupling is optimal for bicausal optimal transport problems with submodular costs and \(\sigma,\bar{\sigma}>0\) provided that, for each equation, either its driving noise has Hurst parameter at least \(1/2\) or its coefficients satisfy
\begin{equation}
	\label{eqn-intro-commutator}
	\partial_x b(t,x)\sigma(t,x)-b(t,x)\partial_x\sigma(t,x)-\partial_t\sigma(t,x)\geq0,
\end{equation}
with the analogous condition for \((\bar b,\bar\sigma)\).
In the proof, we establish that the solution map of the discrete approximation process \(\hat{X}^{n}\) is a \emph{pointwise} increasing functional of the driving Brownian increments \(\Delta B\).
Together with the submodular cost and positivity of the diffusion coefficient, we show that the synchronous coupling is optimal for the bicausal transport problem \(\mathcal{T}_{\bc}(\hat{X}^{n},\hat{Y}^{n}; \hat{c}_{n})\) for any \(n\).
Then the optimality of the synchronous coupling follows from the stability result \eqref{eqn-intro-approx}.
To the best of our knowledge, this is the first general bicausal optimality theorem for scalar fractional SDEs including the rough regime \(H<1/2\).

\subsection{Related literature}
Causal optimal transport was introduced in the pioneering work of \citet{lassalle18Causal} and subsequently developed in \citet{backhoff-veraguas17Causal, acciaio20Causal, bartl24Wasserstein, bartl25Wasserstein, cont2024causal}, among others.
A seminal result of \citet{backhoff-veraguas20All} shows that, for discrete-time processes, the adapted Wasserstein distance metrizes several notions of adapted weak convergence, including Hoover's adapted distribution \citep{hoover84Adapteda}, Aldous's generalized weak topology \citep{aldous81Weak}, Hellwig's information topology \citep{hellwig96Sequential}, and the nested distance \citep{pflug12Distance}.
Causal optimal transport has also been extended to richer causal structures, including causal graphs; see \citet{cheridito25optimala, obloj26graph}.
Applications of causal optimal transport include machine learning \citep{xu20cotgan}, mathematical finance \citep{backhoff-veraguas20Adapted}, and distributionally robust optimization \citep{jiang24Duality, jiang24Sensitivity}.

Our stability theorem is closest to the discrete-time stability results of \citet{eckstein24Computational}.
As pointed out in the introduction, however, the discrete-time theory requires convergence of the marginal processes in the adapted Wasserstein distance.
A direct continuous-time analogue of this condition would exclude basic discretizations such as polygonal interpolations.
Our notion of convergence in adapted Wasserstein distance along  common vanishing time shifts allows one to overcome this constraint.
Moreover, \eqref{eqn-intro-stability} allows the filtrations and costs to vary and uses vanishing time shifts to handle the discontinuous information flow of discrete approximations.
A related relaxation of the adapted Wasserstein distance for continuous-time processes is studied in \citet{bartl25Wasserstein}, where the authors considered
\begin{equation*}
	\widetilde{\AW}(X,Y):=\inf_{\varepsilon>0}\Bigl\{\varepsilon+ \inf_{\pi\in \Pi_{\bc,\varepsilon}(\mathbf{F},\mathbf{G})}\E_{\pi}[d(X,Y)]\Bigr\}
\end{equation*}
and \(\Pi_{\bc,\varepsilon}(\mathbf{F},\mathbf{G})\) denotes the set of \(\varepsilon\)-relaxed bicausal couplings.
It is shown that \(\widetilde{\AW}\) metrizes the continuous-time adapted weak topology.
We remark that if \(\pi\in \Pi_{\bc}(\mathbf{F}\circ \tau, \mathbf{G}\circ \tau)\), then \(\pi\in \Pi_{\bc,\|\tau - \Id\|_{\infty}}(\mathbf{F},\mathbf{G})\).
Thus, convergence in adapted Wasserstein along  vanishing time shifts \(\bt\) implies the convergence in the adapted weak topology.

Optimality of the synchronous coupling has been studied for  stochastic differential equations  in \citet{bion-nadal19Wassersteintype} via stochastic control approach and in \citet{backhoff-veraguas22Adapted, robinson24Bicausal} via time discretization.
For fractional noises, \citet{jiang25transfer} first established synchronous optimality for scalar fractional SDEs in the Young regime \(H>1/2\) with a quadratic cost.
We stress that the existing discretization-based results considered It\^o SDEs and the  discretization  required a tailored Euler schemes to  preserve the  monotonicity of the solution map.
In contrast, we consider Stratonovich SDEs, and more generally, Doss--Sussmann solutions to fractional SDEs.
In this setting, the standard Wong--Zakai approximation naturally yields a monotone discretization scheme.
This enables us to extend the synchronous optimality to the rough regime \(H<1/2\) under \eqref{eqn-intro-commutator}, and to general continuous submodular costs on path space.

\subsection{Outline}

The rest of the paper is organized as follows.
In \cref{sec-prelim}, we fix notations and recall causal and bicausal optimal transport.
Section~\ref{sec-fine} establishes compactness, gluing, and shadowing properties of causal couplings.
Section~\ref{sec-stab} proves the stability of causal coupling sets under vanishing time shifts and the general stability Theorem \ref{thm-stab}.
Section~\ref{sec-disc} applies this theorem to dyadic Wong--Zakai approximations and proves Theorem \ref{thm-approx} for Brownian and fractional Brownian drivers.
Section~\ref{sec-optim} establishes the monotonicity of the scalar Doss--Sussmann solution map and proves synchronous optimality for fractional SDEs with submodular path costs.

\section{Preliminaries}
\label{sec-prelim}

\subsection{Notations}
Let \((\mathcal{M},d_{\mathcal{M}})\) be a complete separable metric space.
We may without loss of generality assume \(d_{\mathcal{M}}\) is a bounded metric.
We denote the space of Borel probability measures on \(\mathcal{M}\) by \(\mathscr{P}(\mathcal{M})\) equipped with the weak topology.
By \(C([0,T]; \mathcal{M})\) we denote the continuous path space on \(\mathcal{M}\) equipped with the uniform metric.
For any sequence of subsets \(\{A_{n}\}_{n\geq 1}\) of \(\mathcal{M}\), the Kuratowski limit inferior is defined as
\begin{equation*}
	\liminf_{n\to \infty} A_{n}= \{x\in \mathcal{M} : \limsup_{n\to\infty} d_{\mathcal{M}}(x,A_{n})=0\};
\end{equation*}
the Kuratowski limit superior is defined as
\begin{equation*}
	\limsup_{n\to \infty} A_{n}= \{x\in \mathcal{M} : \liminf_{n\to\infty} d_{\mathcal{M}}(x,A_{n})=0\}.
\end{equation*}
If the Kuratowski limits inferior and superior agree, we call their common value the Kuratowski limit of the sequence and denote it by \(\lim_{n\to\infty}A_{n}\).

For paths \(x,x'\in C([0,T]; \mathbb{R})\), we define a partial order \(x\preceq x'\) if for any \(t\in[0,T]\) it holds that \(x(t)\leq  x'(t)\).
We say a functional \(\Phi\) on \(C([0,T]; \mathbb{R})\) is \emph{pointwise} increasing if
\begin{equation*}
	\Phi(x)\leq \Phi(x')  \text{ for any } x\preceq x'.
\end{equation*}
Note that we use increasing always in a relaxed sense, i.e., non-decreasing.

Let \(\mathcal{X}\) and \(\mathcal{Y}\) be two complete separable metric spaces.
Given \(\mu\in \mathscr{P}(\mathcal{X})\) and \(\nu\in \mathscr{P}( \mathcal{Y})\),  the set of couplings between \(\mu\) and \(\nu\) is defined as
\begin{equation*}
	\Pi(\mu,\nu):=\{\pi\in \mathscr{P}( \mathcal{X}\times \mathcal{Y}):\pi(\cdot\times \mathcal{Y})=\mu(\cdot)\text{ and } \pi(\mathcal{X}\times\cdot)=\nu(\cdot)\}.
\end{equation*}
For any \(\pi\in \mathscr{P}(\mathcal{X}\times \mathcal{Y})\), we write
\begin{equation*}
	\pi(\md x,\md y)=\pi(\md x) \theta(x,\md y),
\end{equation*}
where  \(\theta\) is the Borel regular disintegration kernel.
The Wasserstein distance on \(\mathscr{P}(\mathcal{X})\) is given by
\begin{equation*}
	\W(\mu,\nu)=\inf_{\pi\in \Pi(\mu,\nu)}\E_{\pi}[d_{\mathcal{X}}(X,Y)].
\end{equation*}

Throughout the paper, we fix a standard Borel probability space \((\Omega,\mathcal{F},\P)\).
We denote the distribution of a random variable \(X\) by \(\L(X)\).
For random variables \(X^{n}\) taking values in \(\mathcal{X}\), we say \(X^{n}\) converges in distribution to \(X\), \(X^{n}\xrightarrow{\D}X\), if and only if
\begin{equation*}
	\L(X^{n})\to \L(X)
\end{equation*}
in the weak topology.

We may consider stochastic processes \(X=(X_{t})_{t\in I}\) on \(\Omega\) in both discrete and continuous time.
Given a filtration \(\mathbf{F}=(\mathcal{F}_{t})_{t\in I}\) with \(\mathcal{F}_{t}\) sub \(\sigma\)-algebra of \(\mathcal{F}\), we often write $(X, \mathbf{F})$ to refer to an $\mathbf{F}$-adapted stochastic process $X$.
If \(\mathbf{F}\) is the raw natural filtration induced by \(X\), then we abbreviate \((X,\mathbf{F})=X\).

\subsection{Causal optimal transport}
Let the time index \(I\) be either a finite set \(\{0,1,\dots,T\}\) in discrete time or \([0,T]\) in continuous time.
\begin{defn}[Causal coupling]
	We say a coupling \(\pi\in\Pi(\P,\P)\) is \emph{causal} from \(\mathbf{F}\) to \(\mathbf{G}\) if for all $t\in I$
	\[\mathcal{F}_{T}\otimes \{\emptyset, \Omega\} \text{ and } \{\emptyset, \Omega\}\otimes\mathcal{G}_{t} \text{ are conditionally independent given } \mathcal{F}_{t}\otimes \{\emptyset, \Omega\} \text{ under } \pi.\]
	We say \(\pi\in\Pi(\P,\P)\) is bicausal if it is causal and \([(x_{1},x_{2})\mapsto (x_{2},x_{1})]_{\#}\pi\) is causal from \(\mathbf{G}\) to \(\mathbf{F}\).
	We denote the set of causal (bicausal) couplings from \(\mathbf{F}\) to \(\mathbf{G}\) by \(\Pi_{\c}(\mathbf{F},\mathbf{G})\) \((\Pi_{\bc}(\mathbf{F},\mathbf{G}))\).
\end{defn}

We notice that the causality condition here only depends on the filtration of the underlying probability space and can be easily extended to the case where the source and target probability spaces \((\Omega,\mathcal{F},P)\) are different.

Let \((X,\mathbf{F})\) and \((Y,\mathbf{G})\) be two adapted stochastic processes on \((\Omega,\mathcal{F},\P)\) with paths in \(\mathcal{X}\) and \(\mathcal{Y}\) respectively.
Let \(c:\mathcal{X}\times \mathcal{Y}\to \mathbb{R}\) be a Borel measurable cost function.
The associated (bi)causal optimal transport problem is given by
\begin{equation*}
	\mathcal{T}_{\bullet}((X,\mathbf{F}),(Y,\mathbf{G}); c):=\inf_{\pi\in \Pi_{\bullet}(\mathbf{F},\mathbf{G})}\E_{\pi}[c(X,Y)],
\end{equation*}
for \(\bullet\in \{\c,\bc\}\).
In particular, if \(\mathcal{X}=\mathcal{Y}\) we write
\begin{equation*}
	\AW((X,\mathbf{F}),(Y,\mathbf{G})):= \mathcal{T}_{\bc}((X,\mathbf{F}),(Y,\mathbf{G}); d_{\mathcal{X}})
\end{equation*}
as the adapted Wasserstein distance between \((X,\mathbf{F})\) and \((Y,\mathbf{G})\).
In discrete time, for processes with their natural filtrations, it is shown that \(\AW\) metrizes the adapted weak topology in  \citet{backhoff-veraguas20All}.
For general filtered processes, we refer the reader to \citet{bartl24Wasserstein}.

\subsection{Doss--Sussmann solution}
\label{subsection:dosssussmannsoln}
In \citet{doss77Liens} and \citet{sussmann78Gap}, a  new solution concept has  been introduced for scalar stochastic differential equation beyond semi-martingale driving noises.
It interprets the solution in a pathwise sense through the so-called Doss--Sussmann transform.
Let \(Z\) be a continuous random noise with \(Z_{0}=0\).
We consider the following equation
\begin{equation}
	\label{eqn-ds}
	\md X_{t}= b(t,X_{t})\md t + \sigma(t,X_{t}) \circ\md Z_{t} \qquad X_{0}=x_{0},
\end{equation}
where \(b\in C^{1}_{b}\) and \(\sigma \in C_{b}^{2}\).
We define the flow \(F:[0,T]\times \mathbb{R}\times \mathbb{R} \to \mathbb{R}\) as the solution to the ordinary differential equation:
\begin{equation*}
	\partial_{z}F(t, y,z):= \sigma(t, F(t,y,z)) \text{ with } F(t,y,0)=y,
\end{equation*}
and further define \(G(t,y,z):= \frac{b(t,F(t,y,z))-\partial_{t}F(t,y,z)}{\partial_{y}F(t,y,z)}\).
The Doss--Sussmann solution is given by
\begin{equation*}
	X_{t}= F(t, Y_{t},Z_{t}),
\end{equation*}
where  \(Y\) solves
\begin{equation*}
	\md Y_{t}=G(t,Y_{t},Z_{t})\md t \text{ with } Y_{0}=x_{0}.
\end{equation*}
We point out that by \citet[Lemma 4.]{doss77Liens} the pathwise solution map is continuous with respect to the uniform convergence of the driver.
For smooth driving noise, it agrees with the classical ODE solution.

\section{Fine properties of causal couplings}
\label{sec-fine}
In this section, we discuss the compactness, gluing, and shadowing for (bi)causal couplings on a standard Borel base space with general filtrations.
Some of these results have been shown in the classical framework, where the probability space is the canonical path space and the filtration is the natural filtration.
For completeness, we present a self-contained proof and adapt to the current framework.

We recall the following from \citet[Lemma 2.5]{beiglbock20denseness}.
\begin{lem}
	\label{lem-topoborelsetdependence}
	Suppose $\pi, \pi_n \in \Pi(P,P)$ for $n \geq 0$. Then $\pi_n$ converges weakly to $\pi$ if and only if for all bounded \(\mathcal{F}\otimes \{\emptyset,\Omega\}\)-random variable \(\varphi\) and bounded \(\{\emptyset,\Omega\}\otimes \mathcal{F}\)-random variable \(\psi\),
	\begin{equation*}
		\lim_{n\to\infty}\E_{\pi_{n}}[\varphi \psi]=\E_{\pi}[\varphi \psi].
	\end{equation*}
	In particular, the topology of $\Pi(P,P)$ depends only on the Borel \(\sigma\)-algebra \(\mathcal{F}\) of $\Omega$, not on the topologies.
\end{lem}

\begin{prop}[Compactness]
	\label{prop-compact}
	The coupling set \(\Pi_{\bullet}(\mathbf{F},\mathbf{G})\) is compact under the weak topology for \(\bullet\in\{\c, \bc\}\).
\end{prop}
\begin{proof}
	We prove the proposition for $\bullet = \c$.
	The case for $\bullet = \bc$ follows in an analogous manner.
	Following from the compactness of \(\Pi(\P,\P)\) and
	\begin{equation*}
		\Pi_{\c}(\mathbf{F},\mathbf{G}) = \bigcap_{t\in I}\Pi_{\c,t}(\mathbf{F},\mathbf{G}):=\bigcap_{t\in I} \Pi_{\c}((\mathcal{F}_{t},\mathcal{F}_{T}), (\mathcal{G}_{t},\mathcal{G}_{T})),
	\end{equation*}
	it suffices to show that each \(\Pi_{\c,t}(\mathbf{F},\mathbf{G})\) is closed.
	We note that $\pi \in \Pi_{\c,t}(\mathbf{F},\mathbf{G})$ if and only if \[
		E_\pi[\Phi_{T} | \mathcal{F}_t \otimes \{ \emptyset , \Omega \}] = E_\pi[\Phi_{T}| \mathcal{F}_t \otimes \mathcal{G}_t],
	\] for all bounded $\mathcal{F}_T\otimes  \{ \emptyset , \Omega \}$-random variable \(\Phi_{T}\).
	This is equivalent to
	\[
		\E_{\pi}[\E_{\pi}[\Phi_{T}|\mathcal{F}_t\otimes  \{ \emptyset , \Omega \} ] \Phi_{t} \Psi_{t}]= \E_{\pi}[\Phi_{T}\Phi_{t}\Psi_{t}]
	\] for any bounded $\mathcal{F}_t\otimes  \{ \emptyset  ,\Omega \}$-random variable \(\Phi_{t}\) and bounded $  \{ \emptyset , \Omega \}\otimes\mathcal{G}_t$-random variable \(\Psi_{t}\).
	We now fix a sequence $(\pi_n)_{n \geq 0} \subseteq \Pi_{\c,t}(\mathbf{F}, \mathbf{G})$ converging weakly to some $\pi \in \Pi(\P,\P)$.
	Since \(\pi_{n}\) is causal, it holds that
	\begin{equation}
		\label{eqn-cau}
		\E_{\pi_{n}}[\E_{\pi_{n}}[\Phi_{T}|\mathcal{F}_t\otimes  \{ \emptyset , \Omega \} ] \Phi_{t} \Psi_{t}]= \E_{\pi_{n}}[\Phi_{T}\Phi_{t}\Psi_{t}].
	\end{equation}
	By choosing $\varphi = \Phi_{T}\Phi_{t}$ and $\psi = \Psi_{t}$, we obtain by Lemma \ref{lem-topoborelsetdependence}
	\[
		\lim_{n\to\infty}\E_{\pi_{n}}[\Phi_{T}\Phi_{t}\Psi_{t}]=\E_{\pi}[\Phi_{T}\Phi_{t}\Psi_{t}].
	\]
	On the other hand, since \(\pi_{n}\) has the same first marginal, there exists a bounded \(\mathcal{F}_{T}\otimes \{\emptyset,\Omega\}\)-random variable \(\hat{\varphi}\) such that
	\begin{equation*}
		\hat{\varphi}= \E_{\pi_{n}}[\Phi_{T}|\mathcal{F}_t\otimes  \{ \emptyset , \Omega \} ] \Phi_{t} \qquad \pi_{n}\text{-a.s.}
	\end{equation*}
	By  applying Lemma \ref{lem-topoborelsetdependence} again with \(\hat{\varphi}\) and \(\psi\) we derive that
	\begin{equation*}
		\lim_{n\to\infty}\E_{\pi_{n}}[\E_{\pi_{n}}[\Phi_{T}|\mathcal{F}_t\otimes  \{ \emptyset , \Omega \} ] \Phi_{t} \Psi_{t}]=\E_{\pi}[\E_{\pi}[\Phi_{T}|\mathcal{F}_t\otimes  \{ \emptyset , \Omega \} ] \Phi_{t} \Psi_{t}]
	\end{equation*}
	Together with \eqref{eqn-cau}, we deduce the causality of \(\pi\).
\end{proof}

The preceding proposition gives the existence of an optimizer for the causal/bicausal transport problem without \emph{any} assumptions on the filtrations $\mathbf{F}$ and $\mathbf{G}$ as long as the base space is standard Borel:

\begin{cor}
	\label{cor-attain}
	Let \(c:\mathcal{X}\times \mathcal{Y}\to \mathbb{R}\) be lower semicontinuous and bounded from below.
	The infimum for \(\mathcal{T}_{\bullet}((X,\mathbf{F}),(Y,\mathbf{G}); c)\) is attained for \(\bullet\in \{\c,\bc\}\).
\end{cor}
\begin{proof}
	We recall the definition
	\begin{equation*}
		\mathcal{T}_{\bullet}((X,\mathbf{F}),(Y,\mathbf{G}); c):=\inf_{\pi\in \Pi_{\bullet}(\mathbf{F},\mathbf{G})}\E_{\pi}[c(X,Y)].
	\end{equation*}
	By \citet[Theorem 13.11]{kechris95classical}, there exists a Polish topology compatible with \((\Omega,\mathcal{F})\) such that \(X\) and \(Y\) are continuous.
	Therefore, \(c(X,Y):\Omega\times \Omega\to \mathbb{R}\)  is lower semi-continuous and bounded from below.
	Together with Proposition \ref{prop-compact} and the Portmanteau theorem we deduce the existence of the optimal coupling.
\end{proof}

Since we consider general filtrations \(\mathbf{F}\) and \(\mathbf{G}\), \((\Omega,\mathcal{F}_{T})\) and \((\Omega,\mathcal{G}_{T})\) need not be standard Borel spaces.
Nevertheless, the following lemma allows us to construct  regular disintegration kernels with respect to \(\mathcal{F}_{T}\) and \(\mathcal{G}_{T}\).

\begin{lem}
	\label{lem-dis}
	Let \(\hat{\pi}\in \mathscr{P}(\Omega\times\Omega, \mathcal{F}_{T}\otimes \mathcal{G}_{T})\) with marginals  \(\P|_{\mathcal{F}_{T}}\) and \(\P|_{\mathcal{G}_{T}}\).
	There exists a regular disintegration kernel \(\hat{\theta}: (\Omega,\mathcal{F}_{T}) \times \mathcal{G}_{T}\to [0,1]\)  such that  \(\hat{\pi}(\md x,\md y)= \hat{\theta}(x,\md y) \P|_{\mathcal{F}_{T}}(\md x)\).
\end{lem}
\begin{proof}
	We first construct \(\pi\) an extension of \(\hat{\pi}\) to \(\mathcal{F}_{T}\otimes \mathcal{F}\).
	Since \((\Omega,\mathcal{F})\) is a standard Borel space, there exists a regular conditional kernel
	\(K:(\Omega,\mathcal{G}_{T})\times \mathcal{F}\to [0,1]\) such that
	\begin{equation*}
		K(y, U)= \E[\1_{U}|\mathcal{G}_{T}] \text{ for any } U\in \mathcal{F}.
	\end{equation*}
	We define \(\pi(\md x,U)= \int K(y, U)\hat{\pi}(\md x, \md y)\) for any \(U\in \mathcal{F}\).
	It is clear that \(\pi\in \mathscr{P}(\Omega\times \Omega, \mathcal{F}_{T}\otimes \mathcal{F})\) and \(\pi|_{\mathcal{F}_{T}\otimes \mathcal{G}_{T}}=\hat{\pi}\).
	Again, since \((\Omega,\mathcal{F})\) is a standard Borel space, \(\pi\) has a regular disintegration kernel \(\theta:(\Omega,\mathcal{F}_{T})\times \mathcal{F} \to [0,1]\) such that
	\begin{equation*}
		\pi(\md x,\md y)= \theta(x,\md y) P(\md x).
	\end{equation*}
	We take \(\hat{\theta}= \theta|_{(\Omega,\mathcal{F}_{T})\times \mathcal{G}_{T}}\) and complete the proof.
\end{proof}

\begin{lem}[Gluing]
	\label{lem-glue}
	Let \(\mathbf{F}\), \(\mathbf{G}\), \(\mathbf{H}\) be three filtrations and $\bullet \in \{ \c, \bc\}$.
	Assume that \(\pi \in \Pi_{\bullet}(\mathbf{F},\mathbf{G})\) and \(\lambda \in \Pi_{\bullet}(\mathbf{G},\mathbf{H})\).
	Then, there exists some $\zeta(\md x,\md y,\md z) \in \mathscr{P}(\Omega\times \Omega\times \Omega, \mathcal{F}\otimes \mathcal{F} \otimes \mathcal{F})$  such that $\zeta(\md x,\md y)|_{\mathcal{F}_T \otimes \mathcal{G}_T} = \pi(\md x,\md y)|_{\mathcal{F}_T \otimes \mathcal{G}_T}, \zeta(\md y,\md z)|_{\mathcal{G}_T \otimes \mathcal{H}_T} = \lambda(\md y,\md z)|_{\mathcal{G}_T \otimes \mathcal{H}_T}$ and $\zeta(\md x,\md z) \in \Pi_{\bullet}(\mathbf{F},\mathbf{H})$.
\end{lem}

\begin{proof}
	We prove the case where $\bullet = \c$, noting that the case for $\bullet = \bc$ follows in a very similar manner.
	By Lemma \ref{lem-dis}, there exists regular disintegration kernels
	$\theta_\pi : (\Omega, \mathcal{F}_T)\times \mathcal{G}_{T}\to [0,1]$ and $\theta_\lambda : (\Omega, \mathcal{G}_T) \times \mathcal{H}_{T}\to [0,1]$ such that  \[
		\pi|_{\mathcal{F}_{T}\otimes \mathcal{G}_{T}}(\md x, \md y) = \theta_\pi(x,\md y) \P|_{\mathcal{F}_{T}}(\md x) \text{ and } \lambda|_{\mathcal{G}_{T}\otimes \mathcal{H}_{T}}(\md y,\md z) = \theta_\lambda(y,\md z) \P|_{\mathcal{G}_{T}}(\md y).
	\]
	We set
	\begin{equation*}
		\hat{\zeta}(\md x,\md y,\md z) := \theta_{\lambda}(y, \md z)\theta_{\pi}(x, \md y) P|_{\mathcal{F}_{T}}(\md x)\in \mathscr{P}(\Omega\times \Omega\times \Omega, \mathcal{F}_{T}\otimes \mathcal{G}_{T}\otimes \mathcal{H}_{T}).
	\end{equation*}
	Let $\Phi_{t}, \Psi_{t}, \Theta_{t}$ denote arbitrary bounded $\mathcal{F}_t$, $\mathcal{G}_t$, $\mathcal{H}_t$-random variables respectively for $t \in I$.
	We note that from our construction
	\begin{equation*}
		\E_{\hat{\zeta}}[\Phi_{T}\Psi_{T}\Theta_{T}]=\E_{\hat{\zeta}}[\Phi_{T}\Psi_{T}\E_{\hat{\zeta}}[\Theta_{T}|\{\emptyset,\Omega\} \otimes \mathcal{G}_T \otimes \{\emptyset,\Omega\}]].
	\end{equation*}
	Thus, we deduce
	\begin{equation*}
		\E_{\hat{\zeta}}[\Theta_{T}|\mathcal{F}_T \otimes \mathcal{G}_T \otimes \{\emptyset,\Omega\}]=\E_{\hat{\zeta}}[\Theta_{T}| \{\emptyset,\Omega\} \otimes \mathcal{G}_T \otimes \{\emptyset,\Omega\}].
	\end{equation*}
	In order to show that under $\hat{\zeta}$, $\mathcal{F}_T \otimes \{\emptyset, \Omega\} \otimes \{ \emptyset, \Omega \}$ and $\{\emptyset, \Omega\}  \otimes \{ \emptyset, \Omega \} \otimes \mathcal{H}_t$ is conditionally independent, given $\mathcal{F}_t \otimes \{\emptyset, \Omega\} \otimes \{ \emptyset, \Omega \}$,  it suffices to verify
	\[
		\E_{\hat{\zeta}} [\E_{\hat{\zeta}}[\Phi_{T}|\mathcal{F}_{t}\otimes \{\emptyset, \Omega\} \otimes \{ \emptyset, \Omega \}] \Phi_{t} \Theta_{t}] = \E_{\hat{\zeta}}[\Phi_{T} \Phi_{t} \Theta_{t}].
	\]
	We compute that \begin{align*}
		  & \quad		\E_{\hat{\zeta}} [\E_{\hat{\zeta}}[\Phi_{T}|\mathcal{F}_{t}\otimes \{\emptyset, \Omega\} \otimes \{ \emptyset, \Omega \}]\Phi_{t} \Theta_{t}]
		\\
		= & \E_{\hat{\zeta}} [\E_{\hat{\zeta}}[\Phi_{T}|\mathcal{F}_{t}\otimes \{\emptyset, \Omega\} \otimes \{ \emptyset, \Omega \}]\Phi_{t} \E_{\hat{\zeta}}[\Theta_{t}|\mathcal{F}_{T}\otimes \mathcal{G}_{T}\otimes \{ \emptyset, \Omega \} ]]         \\
		= & \E_{\hat{\zeta}} [\E_{\hat{\zeta}}[\Phi_{T}|\mathcal{F}_{t}\otimes \{\emptyset, \Omega\} \otimes \{ \emptyset, \Omega \}]\Phi_{t} \E_{\hat{\zeta}}[\Theta_{t}|\{ \emptyset, \Omega \}\otimes \mathcal{G}_{T}\otimes \{ \emptyset, \Omega \} ]] \\
		= & \E_{\hat{\zeta}} [\E_{\hat{\zeta}}[\Phi_{T}|\mathcal{F}_{t}\otimes \{\emptyset, \Omega\} \otimes \{ \emptyset, \Omega \}]\Phi_{t} \E_{\hat{\zeta}}[\Theta_{t}|\{ \emptyset, \Omega \}\otimes \mathcal{G}_{t}\otimes \{ \emptyset, \Omega \} ]] \\
		= & \E_{\hat{\zeta}}[\Phi_{T}\Phi_{t} \E_{\hat{\zeta}}[\Theta_{t}|\{ \emptyset, \Omega \}\otimes \mathcal{G}_{t}\otimes \{ \emptyset, \Omega \} ]]                                                                                                 \\
		= & \E_{\hat{\zeta}}[\Phi_{T}\Phi_{t} \E_{\hat{\zeta}}[\Theta_{t}|\{ \emptyset, \Omega \}\otimes \mathcal{G}_{T}\otimes \{ \emptyset, \Omega \} ]]                                                                                                 \\
		= & \E_{\hat{\zeta}}[\Phi_{T}\Phi_{t} \E_{\hat{\zeta}}[\Theta_{t}|\mathcal{F}_{T}\otimes \mathcal{G}_{T}\otimes \{ \emptyset, \Omega \} ]]                                                                                                         \\
		= & \E_{\hat{\zeta}}[\Phi_{T}\Phi_{t}\Theta_{t}].
	\end{align*}
	The first equality follows from the tower law; the second equality follows from the conditional independence between $\mathcal{F}_T \otimes \{\emptyset, \Omega\} \otimes \{ \emptyset, \Omega \}$ and $\{\emptyset, \Omega\}  \otimes \{ \emptyset, \Omega \} \otimes \mathcal{H}_t$ given $\{ \emptyset, \Omega \} \otimes \mathcal{G}_T \otimes \{ \emptyset, \Omega \}$;
	the third equality follows from \(\lambda\in \Pi_{\c}(\mathbf{G},\mathbf{H})\);
	the fourth equality follows from \(\pi\in \Pi_{\c}(\mathbf{F},\mathbf{G})\);
	the fifth equality follows from \(\lambda\in \Pi_{\c}(\mathbf{G},\mathbf{H})\);
	the sixth equality follows from the conditional independence between $\mathcal{F}_T \otimes \{\emptyset, \Omega\} \otimes \{ \emptyset, \Omega \}$ and $\{\emptyset, \Omega\}  \otimes \{ \emptyset, \Omega \} \otimes \mathcal{H}_t$ given $\{ \emptyset, \Omega \} \otimes \mathcal{G}_T \otimes \{ \emptyset, \Omega \}$;
	the last equality follows from the tower law.

	We next extend $\hat{\zeta}$ to $\mathcal{F}\otimes\mathcal{F}\otimes\mathcal{F}$.
	Since $(\Omega,\mathcal{F})$ is a standard Borel space, we take $K_{\mathcal{F}_{T}},K_{\mathcal{G}_{T}},K_{\mathcal{H}_{T}}$ as regular conditional probability kernels of $\P$ given $\mathcal{F}_{T},\mathcal{G}_{T},\mathcal{H}_{T}$, respectively, and define
	\begin{equation*}
		\zeta(\md x',\md y',\md z'):=\int_{\Omega\times \Omega\times \Omega}K_{\mathcal{F}_{T}}(x,\md x')K_{\mathcal{G}_{T}}(y,\md y')K_{\mathcal{H}_{T}}(z,\md z')\hat{\zeta}(\md x,\md y,\md z).
	\end{equation*}
	It is direct to verify that
	\begin{equation*}
		\zeta|_{\mathcal{F}_{T}\otimes\mathcal{G}_{T}\otimes\mathcal{H}_{T}}=\hat{\zeta}.
	\end{equation*}
	In particular,
	\begin{equation*}
		\zeta(\md x,\md y)|_{\mathcal{F}_{T}\otimes\mathcal{G}_{T}}=\pi(\md x,\md y)|_{\mathcal{F}_{T}\otimes\mathcal{G}_{T}},\qquad \zeta(\md y,\md z)|_{\mathcal{G}_{T}\otimes\mathcal{H}_{T}}=\lambda(\md y,\md z)|_{\mathcal{G}_{T}\otimes\mathcal{H}_{T}}.
	\end{equation*}
	Since the causality condition only depends on its restriction to $\mathcal{F}_{T}\otimes\mathcal{H}_{T}$, and this restriction agrees with the $(x,z)$-marginal of $\hat{\zeta}$, we conclude that
	\begin{equation*}
		\zeta(\md x,\md z)\in\Pi_{\c}(\mathbf{F},\mathbf{H}).
	\end{equation*}

\end{proof}

\begin{lem}
	\label{lem-coup}
	Let $(X, \mathbf{F}), (X', \mathbf{F}'), (Y, \mathbf{G})$ and $(Y', \mathbf{G}')$ be adapted stochastic processes and \(\pi\in \Pi_{\bullet}(\mathbf{F},\mathbf{G})\) for \(\bullet\in\{\c,\bc\}\).
	Then there exists \(\pi'\in \Pi_{\bullet}(\mathbf{F}',\mathbf{G}')\) such that
	\begin{equation*}
		\W((X,Y)_{\#}\pi,(X',Y')_{\#}\pi')\leq \AW((X,\mathbf{F}),(X',\mathbf{F}'))+\AW((Y,\mathbf{G}),(Y',\mathbf{G}')).
	\end{equation*}
\end{lem}
\begin{proof}
	Let \(\lambda\in\Pi_{\bc}(\mathbf{F},\mathbf{F}') \) attain \(\AW((X,\mathbf{F}),(X',\mathbf{F}'))\).
	By Lemma \ref{lem-glue}, there exists \(\zeta(\md x ,\md x', \md y)\) such that \(\tilde{\pi}(\md x',\md y):=\zeta(\md x', \md y)\in \Pi_{\bullet}(\mathbf{F}',\mathbf{G})\) and
	\begin{equation*}
		\zeta(\md x, \md x')|_{\mathcal{F}_{T}\otimes \mathcal{F}_{T}'}= \lambda(\md x, \md x')|_{\mathcal{F}_{T}\otimes \mathcal{F}_{T}'}, \quad  \zeta(\md x,\md y)|_{\mathcal{F}_{T}\otimes \mathcal{G}_{T}}=\pi|_{\mathcal{F}_{T}\otimes \mathcal{G}_{T}}(\md x, \md y).
	\end{equation*}
	Note that \(\eta= [(x,x',y)\mapsto ((X(x),Y(y),X'(x'),Y(y)))]_{\#}\zeta\) gives a coupling between \((X,Y)_{\#}\pi\) and \((X',Y)_{\#}\tilde{\pi}\).
	Thus, we derive
	\begin{equation*}
		\W((X,Y)_{\#}\pi,(X',Y)_{\#}\tilde{\pi})\leq \E_{\eta}[d_{\mathcal{X}}(X(x),X'(x'))+d_{\mathcal{Y}}(Y(y),Y(y))] =  \AW((X,\mathbf{F}),(X',\mathbf{F}')).
	\end{equation*}
	Following similar arguments, there exists \(\pi'\in \Pi_{\bullet}(\mathbf{F}',\mathbf{G}')\) such that
	\begin{equation*}
		\W((X',Y)_{\#}\tilde{\pi}, (X',Y')_{\#}\pi')\leq \AW((Y,\mathbf{G}), (Y',\mathbf{G}')).
	\end{equation*}
	The desired estimate then follows from the triangle inequality of the Wasserstein distance.
\end{proof}
The coupling \(\pi'\) constructed above is called the shadowing of \(\pi\).

\section{Stability of causal optimal transport}
\label{sec-stab}

In this section, we establish a general stability result for (bi)causal optimal transport problems.
In particular, not only the law of the underlying process can be varied, but also the associated filtration and the transport cost.

We first introduce the time shifted filtration.
\begin{defn}
	We say \(\tau:[0,T]\to [0,T]\)  is a time shift if \(\tau\) is increasing and \(\tau(t)\geq t\).
	For a filtration \(\mathbf{F}=(\mathcal{F}_{t})_{t\in[0,T]}\), by \(\mathbf{F}\circ \tau\), we denote the time shifted filtration
	\begin{equation*}
		\mathbf{F}\circ \tau = (\mathcal{F}^{\tau}_{t})_{t\in [0,T]}:= (\mathcal{F}_{\tau(t)})_{t\in[0,T]}.
	\end{equation*}
\end{defn}

The following characterizes the stability of time-shifted causal couplings.

\begin{prop}
	\label{prop-coup}
	Let \(\bt=\{\tau^{1},\tau^{2},\dots\}\) be a sequence of time shifts with \(\|\tau^{n}-\Id\|_{\infty}\to 0\).
	Assume that \(\mathbf{F}\) and \(\mathbf{G}\) are right continuous under completion.
	Then, we have for \(\bullet\in \{\c,\bc\}\)
	\begin{equation*}
		\lim_{n\to\infty}\Pi_{\bullet}(\mathbf{F}\circ \tau^{n}, \mathbf{G}\circ \tau^{n})= \Pi_{\bullet}(\mathbf{F}, \mathbf{G}),
	\end{equation*}
	where the limit is understood as the Kuratowski limit under the Wasserstein distance.
\end{prop}

\begin{proof}
	We focus on the proof of the causal case.
	The bicausal case follows from symmetry.
	Since \(\tau^{n}\) is increasing, we note that the image of \(\tau^{n}\) is contained in \([0,T]\).
	From the definition of the causality, we obtain for any \(n\)
	\begin{equation*}
		\Pi_{\c}(\mathbf{F}, \mathbf{G})\subseteq \Pi_{\c}(\mathbf{F}\circ \tau^{n}, \mathbf{G}\circ \tau^{n}).
	\end{equation*}
	This implies that \(\Pi_{\c}(\mathbf{F},\mathbf{G})\subseteq \liminf_{n\to\infty}\Pi_{\c}(\mathbf{F}\circ \tau^{n},\mathbf{G}\circ \tau^{n})\).

	For the other direction, take \(\pi\in \limsup_{n\to\infty}\Pi_{\c}(\mathbf{F}\circ \tau^{n},\mathbf{G}\circ \tau^{n})\).
	By definition, there exists a subsequence \((n_{k})_{k\geq 1}\) and \(\pi_{k}\in \Pi_{\c}(\mathbf{F}\circ \tau^{n_{k}},\mathbf{G}\circ \tau^{n_{k}})\) such that
	\begin{equation*}
		\lim_{k\to\infty} \W(\pi, \pi_{k})=0.
	\end{equation*}
	We fix \(t\in[0,T]\) and \(s>t\).
	Without loss of generality, we may assume \(\tau^{n_{k}}(t)\leq s\).
	By the  causality of \(\pi_{k}\), for any bounded \(\{\emptyset,\Omega\}\otimes\mathcal{G}_{t}\)-random variable \(\Psi\), we have
	\begin{equation*}
		\E_{\pi_{k}}[\Psi|\mathcal{F}_{T}\otimes \{\emptyset,\Omega\}]=\E_{\pi_{k}}[\Psi|\mathcal{F}_{\tau^{n_{k}}(t)}\otimes \{\emptyset,\Omega\}]= \E_{\pi_{k}}[\Psi|\mathcal{F}_{s}\otimes \{\emptyset,\Omega\}].
	\end{equation*}
	Therefore, \(\pi_{k}\in \Pi_{\c}((\mathcal{F}_{s},\mathcal{F}_{T}),(\mathcal{G}_{t},\mathcal{G}_{s}))\).
	By Proposition \ref{prop-compact}, \(\Pi_{\c}((\mathcal{F}_{s},\mathcal{F}_{T}),(\mathcal{G}_{t},\mathcal{G}_{s}))\) is weakly compact.
	Thus, the closedness leads to
	\begin{equation*}
		\E_{\pi}[\Psi|\mathcal{F}_{T}\otimes \{\emptyset,\Omega\}]=\E_{\pi}[\Psi|\mathcal{F}_{s}\otimes \{\emptyset,\Omega\}].	\end{equation*}
	Let \(s\to t+\).
	The right continuity of the completed filtration \((\prescript{\P}{}{}\mathcal{F}_{t})_{t\in[0,T]}\) guarantees that for any bounded \(\{\emptyset,\Omega\}\otimes\mathcal{G}_{t}\)-random variable \(\Psi\)
	\begin{equation*}
		\E_{\pi}[\Psi|\mathcal{F}_{T}\otimes \{\emptyset,\Omega\}]=\E_{\pi}[\Psi|\mathcal{F}_{t}\otimes \{\emptyset,\Omega\}].
	\end{equation*}
	Since \(t\in[0,T]\) is arbitrary, we derive
	\begin{equation*}
		\limsup_{n\to\infty}\Pi_{\c}(\mathbf{F}\circ \tau^{n}, \mathbf{G}\circ \tau^{n})\subseteq \Pi_{\c}(\mathbf{F}, \mathbf{G}).
	\end{equation*}
\end{proof}

\begin{defn}
	Let \((X^{n},\mathbf{F}^{n})\) and \((X,\mathbf{F})\) be adapted stochastic processes on \(\Omega\).
	Assume that \(\bt=\{\tau^{1},\tau^{2},\dots\}\) is a sequence of time shifts with \(\|\tau^{n}-\Id\|_{\infty}\to 0\).
	We say \((X^{n},\mathbf{F}^{n})\) converges to \((X,\mathbf{F})\) in  adapted Wasserstein along \(\bt\) if
	\begin{equation*}
		\lim_{n\to\infty}\AW((X^{n},\mathbf{F}^{n}\circ \tau^{n}), (X, \mathbf{F}\circ \tau^{n}))=0.
	\end{equation*}
\end{defn}

\begin{rmk}
	For a fixed sequence \(\bt\), strict adapted Wasserstein convergence implies convergence  along \(\bt\).
	Moreover, if \(\varepsilon_{n}=\|\tau^{n}-\Id\|_{\infty}\), a coupling bicausal for \(\mathbf{F}^{n}\circ \tau^{n}\) and \(\mathbf{F}\circ \tau^{n}\) is \(\varepsilon_{n}\)-bicausal for \(\mathbf{F}^{n}\) and \(\mathbf{F}\).
	Thus, convergence along \(\bt\) implies the convergence in the relaxed adapted Wasserstein distance in \citet{bartl25Wasserstein}.
\end{rmk}

Our main stability result is the following.
\begin{thm}[Stability]
	\label{thm-stab}
	Let \((X^{n},\mathbf{F}^{n})\), \((X,\mathbf{F})\) be adapted stochastic processes with paths in \(\mathcal{X}\) and \((Y^{n},\mathbf{G}^{n})\), \((Y,\mathbf{G})\) be stochastic processes with paths in \(\mathcal{Y}\).
	Let \(c,c_{n}:\mathcal{X}\times \mathcal{Y}\to \mathbb{R}_{\geq 0}\) be  non-negative continuous cost functions.
	We assume the following:
	\begin{enumerate}[(i)]
		\item \((X^{n},\mathbf{F}^{n})\) and \((Y^{n},\mathbf{G}^{n})\) converge to \((X,\mathbf{F})\) and  \((Y,\mathbf{G})\) in adapted Wasserstein along the same sequence \(\bt\) with \(\|\tau^{n}-\Id\|_{\infty}\to 0\);
		\item  \(\mathbf{F}\) and \(\mathbf{G}\) are right continuous after completion;
		\item \(c_{n}\) converges to \(c\) uniformly on compact sets;
		\item there exist non-negative functions \(a_{n}\in L^{1}(X^{n})\) and \(b_{n} \in L^{1}(Y^{n})\) such that
		      \begin{equation*}
			      c_{n}(x,y)\leq a_{n}(x)+b_{n}(y) \text{ and } \lim_{R\to \infty} \sup_{n}\{\E[a_{n}(X^{n})\1_{a_{n}(X^{n})\geq R}]+\E[b_{n}(Y^{n})\1_{b_{n}(Y^{n})\geq R}]\}=0.
		      \end{equation*}
	\end{enumerate}
	Then, we have for \(\bullet\in \{\c,\bc\}\)
	\begin{equation*}
		\lim_{n\to \infty} \mathcal{T}_{\bullet}((X^{n},\mathbf{F}^{n}\circ \tau^{n}), (Y^{n},\mathbf{G}^{n}\circ \tau^{n}); c_{n})= \mathcal{T}_{\bullet}((X,\mathbf{F}), (Y,\mathbf{G}); c).
	\end{equation*}
\end{thm}

\begin{proof}
	We focus on the convergence of the causal case, where \(\bullet=\c\).
	The bicausal case follows from similar arguments.
	Without loss of generality, we assume that \(X:\Omega\to \mathcal{X}\) and \(Y:\Omega\to \mathcal{Y}\) are continuous, see \citet[Theorem 13.11]{kechris95classical}.

	\emph{Step 1.} We consider the case that \(c_{n}=c\) is fixed, continuous (on $\mathcal{X} \times \mathcal{Y}$), and bounded.
	By Corollary \ref{cor-attain}, there exists \(\pi_{n}\) attaining \(\mathcal{T}_{\c}((X^{n},\mathbf{F}^{n}\circ \tau^{n}), (Y^{n},\mathbf{G}^{n}\circ \tau^{n}); c)\).
	By Lemma \ref{lem-coup}, there exists \(\pi_{n}'\in \Pi_{\c}(\mathbf{F}\circ \tau^{n},\mathbf{G} \circ \tau^{n})\) such that
	\begin{equation*}
		\mathcal{W}((X^{n},Y^{n})_{\#}\pi_{n}, (X,Y)_{\#}\pi_{n}')\leq \AW((X^{n},\mathbf{F}^{n}\circ \tau^{n}),(X, \mathbf{F} \circ \tau^{n}))+ \AW((Y^{n},\mathbf{G}^{n}\circ \tau^{n}),(Y, \mathbf{G} \circ \tau^{n})).
	\end{equation*}
	Given assumption (i), we deduce
	\begin{equation*}
		\liminf_{n\to\infty} \mathcal{T}_{\c}((X^{n},\mathbf{F}^{n}\circ \tau^{n}), (Y^{n},\mathbf{G}^{n}\circ \tau^{n}); c)=
		\liminf_{n\to\infty} \E_{\pi_{n}}[c(X^{n},Y^{n})]=
		\liminf_{n\to\infty} \E_{\pi_{n}'}[c(X,Y)].
	\end{equation*}
	By Proposition \ref{prop-coup}, any converging subsequence of \(\pi'_{n}\) has a limit in \(\Pi_{\c}(\mathbf{F},\mathbf{G})\).
	Together with the continuity of \(X\) and \(Y\), we derive
	\begin{equation*}
		\mathcal{T}_{\c}((X,\mathbf{F}),(Y,\mathbf{G}); c)\leq \liminf_{n\to\infty}\mathcal{T}_{\c}((X^{n},\mathbf{F}^{n}\circ \tau^{n}), (Y^{n},\mathbf{G}^{n}\circ \tau^{n}); c).
	\end{equation*}

	For the reverse direction, by Corollary \ref{cor-attain} we take \(\pi\) attaining \(\mathcal{T}_{\c}((X,\mathbf{F}),(Y,\mathbf{G}); c)\).
	Again by Lemma \ref{lem-coup}, there exists \(\pi_{n}''\in \Pi_{\c}(\mathbf{F}^{n}\circ \tau^{n}, \mathbf{G}^{n}\circ \tau^{n})\) such that
	\begin{equation*}
		\mathcal{W}((X^{n},Y^{n})_{\#}\pi_{n}'', (X,Y)_{\#}\pi)\leq \AW((X^{n},\mathbf{F}^{n}\circ \tau^{n}),(X, \mathbf{F} \circ \tau^{n}))+ \AW((Y^{n},\mathbf{G}^{n}\circ \tau^{n}),(Y, \mathbf{G} \circ \tau^{n})).
	\end{equation*}
	Therefore, we deduce that
	\begin{align*}
		\limsup_{n\to\infty}\mathcal{T}_{\c}((X^{n},\mathbf{F}^{n}\circ \tau^{n}), (Y^{n},\mathbf{G}^{n}\circ \tau^{n}); c)\leq \limsup_{n\to\infty}\E_{\pi_{n}''}[c(X^{n},Y^{n})] & =\E_{\pi}[c(X,Y)]                                    \\
		                                                                                                                                                                           & =\mathcal{T}_{\c}((X,\mathbf{F}),(Y,\mathbf{G}); c).
	\end{align*}

	\emph{Step 2.} Suppose now \(c_{n}\) converges to \(c\) uniformly on compact sets such that \(\|c_{n}\|_{\infty}\leq 1\) for all \(n\geq 0\).
	We note that assumption (i) implies that \((X^{n}_{\#}\P)_{n\geq 1}\) and \((Y^{n}_{\#}\P)_{n\geq 1}\) are tight sequences.
	For any \(\varepsilon>0\), there exist compact sets \(K_{\varepsilon}\subseteq \mathcal{X}\) and \(L_{\varepsilon}\subseteq \mathcal{Y}\) such that
	\begin{equation*}
		\sup_{n\geq 1} \{P(X^{n}\notin K_{\varepsilon}) + P (Y^{n}\notin L_{\varepsilon})\}<\varepsilon.
	\end{equation*}
	For \(n\) sufficiently large, we have  \(|c_{n}-c|\leq \varepsilon\) on \(K_{\varepsilon}\times L_{\varepsilon}\).
	Therefore, we have
	\begin{align*}
		|\E_{\pi}[c_{n}(X^{n},Y^{n})-c(X^{n},Y^{n})]| & \leq \E_{\pi}[|c_{n}(X^{n},Y^{n})-c(X^{n},Y^{n})|\1_{(X^{n},Y^{n})\in K_{\varepsilon}\times L_{\varepsilon}}]                 \\
		                                              & \quad + (\|c_{n}\|_{\infty}+\|c\|_{\infty})	\pi((X^{n},Y^{n})\notin K_{\varepsilon}\times L_{\varepsilon})\leq 3 \varepsilon.
	\end{align*}
	Hence, we obtain
	\begin{equation*}
		\limsup_{n\to\infty}| \mathcal{T}_{\c}((X^{n},\mathbf{F}^{n}\circ \tau^{n}), (Y^{n},\mathbf{G}^{n}\circ \tau^{n}); c_{n}) -\mathcal{T}_{\c}((X^{n},\mathbf{F}^{n}\circ \tau^{n}), (Y^{n},\mathbf{G}^{n}\circ \tau^{n}); c)|\leq 3 \varepsilon.
	\end{equation*}
	By taking \(\varepsilon\to 0\), we show the convergence of the transport cost from \emph{Step 1}.

	\emph{Step 3.} We then consider general unbounded costs \(c_{n}\).
	Since \(c_{n}(x,y)\leq a_{n}(x)+b_{n}(y)\), we have
	\begin{align*}
		0 & \leq \E_{\pi}[c_{n}(X^{n},Y^{n})- (c_{n}\wedge R)(X^{n},Y^{n})]                                       \\
		  & \leq 2\E_{\pi}[\max\{a_{n}(X^{n}),b_{n}(Y^{n})\}\1_{\{2\max\{a_{n}(X^{n}),b_{n}(Y^{n})\}\geq R\}} ]   \\
		  & \leq 2 \E[a_{n}(X^{n})\1_{\{a_{n}(X^{n})\geq R/2\}}]+2 \E[b_{n}(Y^{n})\1_{\{b_{n}(Y^{n})\geq R/2\}}].
	\end{align*}
	By assumption (iv), for any \(\varepsilon>0\), there exists \(R>0\) such that for any coupling \(\pi\in \Pi(\P,\P)\)
	\begin{equation*}
		0\leq  \E_{\pi}[c_{n}(X^{n},Y^{n})- (c_{n}\wedge R)(X^{n},Y^{n})]  \leq \varepsilon.
	\end{equation*}
	In particular, this implies
	\begin{equation*}
		0\leq \mathcal{T}_{\c}((X^{n},\mathbf{F}^{n}\circ \tau^{n}),(Y^{n},\mathbf{G}^{n}\circ \tau^{n}); c_{n})- \mathcal{T}_{\c}((X^{n},\mathbf{F}^{n}\circ \tau^{n}),(Y^{n},\mathbf{G}^{n}\circ \tau^{n}); c_{n}\wedge R)\leq \varepsilon.
	\end{equation*}
	Similar to \emph{Step 1}, by Lemma \ref{lem-coup}, for any \(\pi\in \Pi_{\c}(\mathbf{F},\mathbf{G})\), there exists \(\pi_{n}\in \Pi_{\c}(\mathbf{F}^{n}\circ \tau^{n},\mathbf{G}^{n}\circ \tau^{n})\) such that
	\begin{equation*}
		\mathcal{W}((X^{n},Y^{n})_{\#}\pi_{n}, (X,Y)_{\#}\pi)\to 0.
	\end{equation*}
	By Skorokhod representation theorem, we may assume \(X^{n}\) converges to \(X\) and \(Y^{n}\) converges to $Y$ almost surely.
	This further implies that \(c_{n}(X^{n},Y^{n})\) converges to \(c(X,Y)\) almost surely.
	Together with assumption (iv), we have \((c_{n}(X^{n},Y^{n}))_{n\geq 1}\) is uniformly integrable, and hence
	\begin{equation*}
		\E_{\pi}[c(X,Y)-(c\wedge R)(X,Y)]=\lim_{n\to\infty}\E_{\pi_{n}}[c_{n}(X^{n},Y^{n})-(c_{n}\wedge R)(X^{n},Y^{n})]\leq \varepsilon.
	\end{equation*}
	It follows that
	\begin{equation*}
		0\leq \mathcal{T}_{\c}((X,\mathbf{F}),(Y,\mathbf{G}); c)- \mathcal{T}_{\c}((X,\mathbf{F}),(Y,\mathbf{G}); c\wedge R)\leq \varepsilon.
	\end{equation*}
	Applying \emph{Step 2} to cost \(c_{n}\wedge R\) and sending \(\varepsilon\to0\), we derive
	\begin{equation*}
		\lim_{n\to\infty} \mathcal{T}_{\c}((X^{n},\mathbf{F}^{n}\circ \tau^{n}), (Y^{n},\mathbf{G}^{n}\circ \tau^{n}); c_{n})= \mathcal{T}_{\c}((X,\mathbf{F}), (Y,\mathbf{G}); c).
	\end{equation*}
\end{proof}

\section{Discrete approximation}
\label{sec-disc}
In this section, we demonstrate how to apply Theorem \ref{thm-stab} to approximate a continuous-time (bi)causal optimal transport by a sequence of discrete-time problems.
Our main focus is Doss--Sussmann solutions to scalar stochastic differential equations:
\begin{equation*}
	\left\{
	\begin{aligned}
		\md X_{t} & = b(t,X_{t}) \md t + \sigma(t,X_{t}) \circ \md Z_{t}\text{ with } X_{0}=x_{0},                         \\
		\md Y_{t} & = \bar{b}(t,Y_{t}) \md t + \bar{\sigma}(t,Y_{t}) \circ \md \bar{Z}_{t} \;\; \text{ with } Y_{0}=y_{0},
	\end{aligned}
	\right.
\end{equation*}
where \(Z\) and \(\bar{Z}\) are random continuous noises whose  natural filtrations are right continuous after completion.

Our aim is to approximate  the problem between continuous-time processes \(X\) and \(Y\) to the problem between discrete-time processes.
We introduce our discretization scheme.
Let \(D_{n}=\{t_{k}=kT2^{-n}: 0\leq k\leq 2^{n}\}\) be the dyadic partition of \([0,T]\).
By \(p_{n}:C([0,T]; \mathbb{R}^{d})\to (\mathbb{R}^{d})^{|D_{n}|}\) we denote the finite marginal projection \(x\mapsto (x_{t_{0}},\dots,x_{t_{|D_{n}|-1}})\);
by \(l_{n}:(\mathbb{R}^{d})^{|D_{n}|}\to C([0,T]; \mathbb{R}^{d})\) we denote the linear interpolation map given by \((x_{t_{0}},\dots,x_{t_{|D_{n}|-1}})\mapsto x\) where
\begin{equation*}
	x_{t}:= \frac{t_{k+1}-t}{t_{k+1}-t_{k}}x_{t_{k}}+\frac{t-t_{k}}{t_{k+1}-t_{k}} x_{t_{k+1}} \text{ for }  t\in[t_{k},t_{k+1}] \text{ and } 0\leq k\leq 2^{n} - 1.
\end{equation*}

In the following, we use a \(\hat{}\) to denote discrete-time processes indexed by \(D_{n}\).
Specifically, let \(\hat{Z}^{n}=(\hat{Z}^{n}_{t})_{t\in D_{n}}\) be a discrete-time process which will be specified later.
We define \(Z^{n}\) as the piecewise linear interpolation of \(\hat{Z}^{n}\), given by
\begin{equation*}
	Z^{n}=l_{n}(\hat{Z}^{n}).
\end{equation*}
We define \(X^{n}\) as the solution to
\begin{equation*}
	\md X^{n}_{t}= b(t,X^{n}_{t}) \md t + \sigma(t,X^{n}_{t}) \md Z^{n}_{t} \text{ with } X^{n}_{0}=x_{0},
\end{equation*}
and set \(\hat{X}^{n}= (X^{n}_{t})_{t\in D_{n}}\).
Similarly, we construct \(\bar{Z}^{n}\), \(Y^{n}\), and \(\hat{Y}^{n}\).

By \(\mathbf{F},\mathbf{F}^{n}, \hat{ \mathbf{F}}^{n}\) we denote the natural filtration of \(X,X^{n},\hat{X}^{n}\), and by \(\mathbf{G},\mathbf{G}^{n},\hat{\mathbf{G}}^{n}\) we denote the natural filtration of \(Y,Y^{n},\hat{Y}^{n}\).
We define a sequence of time shifts \(\bd=\{\delta_{1},\delta_{2},\dots\}\) where \(\delta_{n}:[0,T]\to [0,T]\) is given by
\begin{equation*}
	\delta_{n}(0)=0\text{ and }\delta_{n}(t)= t_{k+1} \text{ for any } t\in(t_{k},t_{k+1}].
\end{equation*}

\begin{asmp}[Coefficients and cost]
	\label{asmp-wp}
	We assume the following:
	\begin{enumerate}[(i)]
		\item \(\sigma, \bar{\sigma}>0.\)
		\item  \(b, \bar{b}\in C^{1}_{b}([0,T]\times \mathbb{R})\) and \(\sigma, \bar{\sigma}\in C^{2}_{b}([0,T]\times \mathbb{R})\).
		\item \(c:C([0,T]; \mathbb{R})\times C([0,T]; \mathbb{R})\to \mathbb{R}\) is continuous, non-negative, and for some \(K,\alpha>0\)
		      \[
			      |c(x,y)|\leq K(1+\|x\|_{\infty}^{\alpha}+\|y\|_{\infty}^{\alpha}).
		      \]
	\end{enumerate}
\end{asmp}

\begin{asmp}[Marginal convergence]
	\label{asmp-conv}
	We assume the following:
	\begin{enumerate}[(i)]
		\item The synchronous coupling \((\Id,\Id)_{\#}\P\) is a bicausal coupling in both \(\Pi_{\bc}(\mathbf{F}^{n}\circ \delta_{n}, \mathbf{F}\circ \delta_{n})\) and \(\Pi_{\bc}(\mathbf{G}^{n}\circ \delta_{n},\mathbf{G}\circ \delta_{n})\).
		\item \(\lim_{n\to\infty}\E[\|X-X^{n}\|_{\infty}\wedge 1]= \lim_{n\to\infty}\E[\|Y-Y^{n}\|_{\infty}\wedge 1]=0\).
		\item Uniform moment estimate: \(\sup_{n}\E\Bigl[\|X^{n}\|_{\infty}^{\beta}+\|Y^{n}\|_{\infty}^{\beta}\Bigr]<\infty\) for some \(\beta>\alpha\).
	\end{enumerate}
\end{asmp}

\begin{thm}
	\label{thm-approx}
	Let \(\hat{X}^{n}=(\hat{X}^{n}_{t})_{t\in D_{n}}\) and \(\hat{Y}^{n}=(\hat{Y}^{n}_{t})_{t\in D_{n}}\) be discrete-time stochastic processes.
	We define cost \(\hat{c}_{n}\) as \(\hat{c}_{n}=c(l_{n}, l_{n})\).
	Under Assumptions \ref{asmp-wp} and \ref{asmp-conv}, it holds that for \(\bullet\in \{\c, \bc\}\)
	\begin{equation*}
		\mathcal{T}_{\bullet}(X,Y; c)= \lim_{n\to \infty}\mathcal{T}_{\bullet}(\hat{X}^{n},\hat{Y}^{n}; \hat{c}_{n}).
	\end{equation*}
\end{thm}

\begin{proof}

	\emph{Step 1.} We first show that
	\begin{equation*}
		\mathcal{T}_{\bullet}(\hat{X}^{n},\hat{Y}^{n};\hat{c}_{n})=\mathcal{T}_{\bullet}(X^{n},Y^{n};c_{n}),
	\end{equation*}
	where
	\begin{equation}
		\label{eqn-cn}
		c_{n}(x,y):=c(l_{n}\circ p_{n}(x),l_{n}\circ p_{n}(y)).
	\end{equation}
	We note that by our construction \(\hat{X}^{n}=p_{n}(X^{n})\) and \(\hat{Y}^{n}=p_{n}(Y^{n})\).
	Therefore, it suffices to show that for \(\bullet\in \{\c,\bc\}\)
	\begin{equation*}
		\Pi_{\bullet}(\hat{\mathbf{F}}^{n},\hat{\mathbf{G}}^{n})= \Pi_{\bullet}(\mathbf{F}^{n},\mathbf{G}^{n}).
	\end{equation*}
	Since \(\sigma\) is non-degenerate, one can verify \(\mathbf{F}^{n}= \mathbf{F}^{Z^{n}}\) and \(\hat{\mathbf{F}}^{n}= \mathbf{F}^{\hat{Z}^{n}}\).
	In particular, \(\mathbf{F}^{n}\) is piecewise constant as \(Z^{n}= l_{n}(\hat{Z}^{n})\), and we have \(\mathcal{F}^{n}_{t}= \mathcal{F}^{\hat{Z}^{n}}_{t}=\hat{\mathcal{F}}^{n}_{t}\) for any \(t\in D_{n}\). With an analogous statement for $\mathbf{G}^n$ and $\hat{\mathbf{G}}^n$, the identity follows immediately  from the definition of causal couplings.

	\emph{Step 2.} We show that \(((X^{n}, \mathbf{F}^{n}), (Y^{n},\mathbf{G}^{n}) ,c_{n})\) satisfies the assumptions in Theorem \ref{thm-stab}.
	We verify that \((X^{n},\mathbf{F}^{n})\) converges to \((X,\mathbf{F})\) and \((Y^{n},\mathbf{G}^{n})\) converges to \((Y,\mathbf{G})\) in adapted Wasserstein along the same sequence \(\bd\).
	Recall that the adapted Wasserstein distance is induced by a bounded metric.
	By Assumption \ref{asmp-conv}, we have
	\begin{equation*}
		\limsup_{n\to\infty}\AW((X^{n},\mathbf{F}^{n}\circ \delta_{n}), (X,\mathbf{F}\circ \delta_{n}))\leq \limsup_{n\to\infty}\E[\|X-X^{n}\|_{\infty}\wedge 1]=0.
	\end{equation*}
	Therefore, \((X^{n},\mathbf{F}^{n})\) converges to \((X,\mathbf{F})\) in adapted Wasserstein along the sequence \(\bd\).
	The same also holds for \(Y\).

	We next verify that \(c_{n}\) defined in \eqref{eqn-cn} converges to \(c\) uniformly on compact sets.
	Let \(K\) be a compact subset of \(C([0,T]; \mathbb{R}^{d})\times C([0,T]; \mathbb{R}^{d})\).
	By Arzel\`a--Ascoli theorem, paths \((x,y)\in K\) are uniformly bounded and equicontinuous.
	We set
	\begin{equation*}
		\tilde{K}=K\cup \{(l_{n}\circ p_{n}(x), l_{n}\circ p_{n}(y)): (x,y)\in K , n\geq 1\}.
	\end{equation*}
	It is direct to show that \(\tilde{K}\) is again uniformly bounded and equicontinuous, and hence precompact.
	Therefore, for any \(\varepsilon>0\) there exists \(\varepsilon'\) such that
	\begin{equation*}
		|c(x,y)-c(x',y')|\leq \varepsilon
	\end{equation*}
	for any \((x,y),(x',y')\in \tilde{K}\) with \(\|(x,y)-(x',y')\|_{\infty}\leq \varepsilon'\).
	By equicontinuity on \(K\), there exists \(n\) sufficiently large such that
	\begin{equation*}
		\|(x,y)-(l_{n}\circ p_{n}(x), l_{n}\circ p_{n}(y))\|_{\infty}\leq \varepsilon' \text{ for all } (x,y)\in K.
	\end{equation*}
	Therefore,  the uniform convergence follows from
	\begin{equation*}
		\limsup_{n\to\infty}\sup_{(x,y)\in K}|c(x,y)-c_{n}(x,y)| = \limsup_{n\to\infty}\sup_{(x,y)\in K} |c(x,y)-c(l_{n}\circ p_{n}(x), l_{n}\circ p_{n}(y))| \leq \varepsilon,
	\end{equation*}
	and sending \(\varepsilon\to 0\).

	We take \(a_{n}(x)=C(1+\|x\|_{\infty}^{\alpha})\) and \(b_{n}(y)=C(1+\|y\|_{\infty}^{\alpha})\).
	By definition \eqref{eqn-cn}, we have \(c_{n}(x,y)\leq a_{n}(x)+b_{n}(y)\).
	It follows from Assumption \ref{asmp-conv}(iii) that
	\begin{align*}
		\E[a_{n}(X^{n})\1_{a_{n}(X^{n})\geq C(1+R)}] & \leq C\E[1+\|X^{n}\|_{\infty}^{\beta}] R^{-\frac{\beta-\alpha}{\alpha}}.\end{align*}
	Therefore, we deduce
	\begin{equation*}
		\lim_{R\to \infty}\sup_{n}\{\E[a_{n}(X^{n})\1_{a_{n}(X^{n})\geq R}]+\E[b_{n}(Y^{n})\1_{b_{n}(Y^{n})\geq R}]\}=0.
	\end{equation*}
	By Theorem \ref{thm-stab}, we deduce for \(\bullet\in\{\c,\bc\}\)
	\begin{equation*}
		\lim_{n\to\infty} \mathcal{T}_{\bullet}((X^{n},\mathbf{F}^{n}\circ \delta_{n}), (Y^{n},\mathbf{G}^{n}\circ \delta_{n}); c_{n})= \mathcal{T}_{\bullet}((X,\mathbf{F}), (Y,\mathbf{G}); c).
	\end{equation*}
	Since \(\mathbf{F}^{n}\circ \delta_{n}= \mathbf{F}^{n}\) and \(\mathbf{G}^{n}\circ \delta_{n}= \mathbf{G}^{n}\), we complete the proof.
\end{proof}

\subsection{Brownian driver}
We focus on Brownian driving noise: \(Z=B\) is a standard Brownian motion.
We give an explicit construction of \(X^{n}\) and \(Z^{n}\) and sufficient conditions for Assumption  \ref{asmp-conv} to hold.
Set \(\hat{Z}^{n}\) as \(\hat{Z}^{n}_{t_{k}}=B_{t_{k}}\) for  \(t_{k}\in D_{n}\), and we retrieve the Wong--Zakai approximation \(X^{n}\) given by
\begin{equation*}
	\md X^{n}_{t}= b(t,X^{n}_{t})\md t + \sigma(t,X^{n}_{t}) \frac{B_{t_{k+1}}-B_{t_{k}}}{ t_{k+1}-t_{k}} \md t \text{ for } t\in [t_{k},t_{k+1}].
\end{equation*}
Similarly, we construct \(Y^{n}\).

The following lemma will be useful for verifying the causality conditions.
\begin{lem}
	\label{lem-markov}
	Let \(M\) be a Markov process on \(\Omega\)  and \(M^{n}=l_{n}\circ p_{n}(M)\).
	We denote the natural filtrations of \(M\) and \(M^{n}\) by \(\mathbf{F}\) and \(\mathbf{F}^{n}\) respectively.
	Then the synchronous coupling \(\pi_{0}=(\Id,\Id)_{\#}P\in \Pi_{\bc}(\mathbf{F}\circ \delta_{n},\mathbf{F}^{n}\circ \delta_{n})\).
\end{lem}
\begin{proof}
	It suffices to show that \(\pi_{0}\in  \Pi_{\c}(\mathbf{F}\circ \delta_{n},\mathbf{F}^{n}\circ \delta_{n})\) and \(\pi_{0}\in  \Pi_{\c}(\mathbf{F}^{n}\circ \delta_{n},\mathbf{F}\circ \delta_{n})\).
	The former follows from the definition of causality and the fact that \(\mathbf{F}\circ \delta_{n}\) is an enlargement of \(\mathbf{F}^{n} \circ \delta_{n}\).
	Let \(t\in(t_{k},t_{k+1}]\) and \(g\in C_{b}\) be \(\mathcal{F}_{\delta_{n}(t)}\)-measurable.
	Then, by the Markov property, \(\sigma((M_{t_{i}})_{i>k+1})\) is conditionally independent of \(\mathcal{F}_{\delta_{n}(t)}=\mathcal{F}_{t_{k+1}}\) given \(\sigma((M_{t_{i}})_{0\leq i \leq k+1})\).
	We denote the copy of \(M\) on the product space \(\Omega\times \Omega\) by \(M'\).
	This yields
	\begin{equation*}
		\E_{\pi_{0}}[g(M')|\mathcal{F}^{n}_{\delta_{n}(T)} \times \{ \emptyset, \Omega \} ]=\E_{P}[g(M)|(M_{t_{i}})_{i\in D_{n}}]=\E_{P}[g(M)| (M_{t_{i}})_{0\leq i\leq k+1}]\in \mathcal{F}^{n}_{t_{k+1}}= \mathcal{F}^{n}_{\delta_{n}(t)}.
	\end{equation*}
	Therefore, we have shown that \(\pi_{0}\in \Pi_{\c}(\mathbf{F}^{n}\circ \delta_{n},\mathbf{F}\circ \delta_{n})\).
\end{proof}

\begin{prop}
	\label{prop-brownian}
	Under Assumption \ref{asmp-wp}, \((X^{n},Y^{n})\) constructed above satisfies Assumption \ref{asmp-conv}.
\end{prop}
\begin{proof}
	We verify the conditions in Assumption \ref{asmp-conv} one by one.

	For (i), we apply Lemma \ref{lem-markov} by taking \(M\) as the driving Brownian motion \(B\).
	Since \(\sigma\) is non-degenerate, the natural filtration of the strong solution \(X\) agrees with the natural filtration of the driving noise \(B\), and the natural filtration of \(X^{n}\) agrees with the natural filtration of \(l_{n}\circ p_{n}(B)\).
	This yields the required bicausal couplings.

	For (ii), from the classical result of Wong--Zakai approximation, we note that \(X^{n}\) converges to \(X\) almost surely.
	By the dominated convergence theorem, this shows the convergence in expectation.

	For (iii),  since \(b\in C^{1}_{b}\) and \(\sigma\in C^{2}_{b}\), the uniform estimates on the moment follow from \citet[Lemma 7.2]{ikeda14stochastic}.
\end{proof}

We remark that the same Wong--Zakai approximation also works with multi-dimensional Stratonovich SDEs.
This implies that the approximating solution \(X^{n}\) converges to the Stratonovich SDE \(X\) in the uniform norm almost surely.
Therefore, Assumption \ref{asmp-conv} is still satisfied, and we can approximate original continuous-time problem by a sequence of discrete ones using Theorem \ref{thm-approx}.

\subsection{Fractional Brownian driver}
For general fractional Brownian noise, we focus on the scalar case \(d=1\).
We take the Molchan--Golosov representation
\begin{equation*}
	Z_{t}=B^{H}_{t}= \int_{0}^{t}k_{H}(t,s)\md B_{s},
\end{equation*}
where \(k_{H}\) is the Molchan--Golosov kernel \citep{molchan69Gaussian, decreusefond99Stochastic} given by
\begin{equation}
	\label{eqn-Molchan-Golosov}
	k_{H}(t,s)= c_H (t-s)^{H-1/2} F(H-1/2,1/2-H,H+1/2,1-t/s)\1_{\{s\leq t\}},
\end{equation} where
\(F(a,b,c,z)\) is the Gaussian hypergeometric function \(
F(a,b,c,z)=\sum_{n=0}^{\infty}\frac{(a)_{n}(b)_{n}}{(c)_{n}}\frac{z^{n}}{n!}\), \((x)_{n}=\Gamma(x+n)/\Gamma(x)\) is the Pochhammer symbol and $c_H > 0$ is some normalizing constant.
In particular, \(\mathbf{F}^{B^{H}}= \mathbf{F}^{B}\).
For any \(t_{k}\in D_{n}\), we set
\begin{equation}
	\label{eqn-zhat}
	\hat{Z}^{n}_{t_{k}}=\frac{2^{n}}{T}\sum_{i=0}^{k-1}\int_{t_{i}}^{t_{i+1}}k_{H}(t_{k},s)\md s (B_{t_{i+1}}-B_{t_{i}}).
\end{equation}
With this choice, we derive a discrete approximation
\begin{equation*}
	\md X^{n}_{t}= b(t,X^{n}_{t}) \md t + \sigma(t,X^{n}_{t}) \md Z^{n}_{t} \text{ with } X^{n}_{0}=x_{0},
\end{equation*}
where
\begin{equation*}
	Z^{n}_{t}:= l_n(\hat{Z}^n)_t = \frac{t_{k+1}-t}{t_{k+1}-t_{k}}\hat{Z}^{n}_{t_{k}}+\frac{t-t_{k}}{t_{k+1}-t_{k}} \hat{Z}^{n}_{t_{k+1}} \text{ for }  t\in[t_{k},t_{k+1}], t_{k},t_{k+1}\in D_{n}.
\end{equation*}

\begin{prop}
	\label{prop-fractional}
	Under Assumption \ref{asmp-wp}, \((X^{n},Y^{n})\) constructed above satisfies Assumption \ref{asmp-conv}.
\end{prop}
\begin{proof}
	We verify the conditions in Assumption \ref{asmp-conv} below.

	For (i), we apply Lemma \ref{lem-markov} again.
	Since \(Z^{n}\) is piecewise linear and \(\sigma\) is non-degenerate, it is direct to show \(\mathbf{F}^{n}=\mathbf{F}^{Z^{n}}\).
	Moreover, from our construction \eqref{eqn-zhat}, this filtration is the same as the natural filtration of \(l_{n}\circ p_{n}(B)\).
	For \(\mathbf{F}\), we note the chain rule for the Doss--Sussmann solution leads to
	\begin{equation*}
		f(t,X_{t})-f(0,x_{0})= \int_{0}^{t} [\partial_{s}f(s,X_{s}) +\partial_{x}f(s,X_{s}) b(s,X_{s}) ]\md s  + \int_{0}^{t} \partial_{x}f(s,X_{s}) \sigma(s,X_{s})\circ \md Z_{s}.
	\end{equation*}
	Taking \(f(t,x)= \int_{0}^{x}\frac{1}{\sigma(t, \xi)}\md \xi\) yields
	\begin{equation*}
		Z_{t}=f(t,X_{t})-f(0,x_{0})- \int_{0}^{t} [\partial_{s}f(s,X_{s}) +\partial_{x}f(s,X_{s}) b(s,X_{s}) ]\md s \in \mathcal{F}_{t}.
	\end{equation*}
	Therefore, we derive \(\mathbf{F}=\mathbf{F}^{Z}\).
	Since we take the canonical representation of the fractional Brownian motion, this filtration is the same as the natural filtration of the driving Brownian motion \(B\).
	By applying Lemma \ref{lem-markov} with \(M=B\), we verify (i).

	For (ii),  it suffices to show that \(\|Z^{n}-Z\|_{\infty}\to 0\) almost surely, since the Doss--Sussmann solution map is continuous in the uniform topology.
	We note that the fractional Brownian motion as a \(C([0,T]; \mathbb{R})\)-valued map \(Z:\Omega\to C([0,T]; \mathbb{R})\) is strongly measurable.
	We define \(\tilde{Z}^{n}:=\E[Z| \mathcal{F}^{n}_{T}]\) as the \(C([0,T]; \mathbb{R})\)-valued conditional expectation.
	A key observation is that
	\begin{equation*}
		Z^{n}_{t_{k}}= \tilde{Z}_{t_{k}}^{n}=\E[Z_{t_{k}}| \mathcal{F}^{n}_{t_{k}}].
	\end{equation*}
	Since \(Z^{n}= l_{n}(\hat{Z}^{n})\), we have
	\begin{align*}
		\|Z^{n}-Z\|_{\infty} & \leq \| l_{n}(\hat{Z}^{n}) - l_{n}\circ p_{n}(Z)\|_{\infty}+ \|l_{n}\circ p_{n}(Z)-Z\|_{\infty}            \\
		                     & = \| l_{n}\circ p_{n}( \tilde{Z}^{n})- l_{n}\circ p_{n} (Z)\|_{\infty}+ \|l_{n}\circ p_{n}(Z)-Z\|_{\infty} \\
		                     & \leq \| \tilde{Z}^{n}-Z\|_{\infty}+\|l_{n}\circ p_{n}(Z)-Z\|_{\infty}.
	\end{align*}
	Since \(Z\in L^{1}(\Omega, C([0,T]; \mathbb{R}))\), and \(\mathcal{F}^{n}_{T}\) increases to \(\mathcal{F}_{T}\), the Banach-valued Martingale convergence theorem (see e.g., \cite[Theorem 11.21]{van2008stochastic}) implies that as \(n\to\infty\) the first term converges to 0 almost surely.
	For the second term, since  \(Z\) is continuous, it also vanishes as \(n\to\infty\).

	For (iii), 	in the scalar case, we recall the Doss--Sussmann solution introduced in Section \ref{subsection:dosssussmannsoln} given by
	\begin{equation*}
		X^{n}_{t}= F(t, Y^{n}_{t},Z^{n}_{t}),
	\end{equation*}
	where  \(Y^{n}\) solves
	\begin{equation*}
		\md Y^{n}_{t}=G(t,Y^{n}_{t},Z^{n}_{t})\md t \text{ with } Y^{n}_{0}=x_{0},
	\end{equation*} and \(F,G:[0,T]\times \mathbb{R}\times \mathbb{R} \to \mathbb{R}\) are as defined in the aforementioned section.
	Since \(\partial_{x}\sigma\) is bounded, we obtain \(\partial_{y}F \geq  \exp(-C|z|)\).
	Taking the derivative \(\partial_{z}(\frac{\partial_{t}F}{\partial_{y}F})\), we obtain
	\begin{equation*}
		\Bigg|\partial_{z}\Bigl(\frac{\partial_{t}F}{\partial_{y}F}\Bigr)\Bigg|=\Bigg|\frac{\partial_{t}\sigma(t,F)}{\partial_{y}F}\Bigg|\leq C \exp(C|z|),
	\end{equation*}
	which yields \( \Bigl|\frac{\partial_{t}F}{\partial_{y}F}\Bigr|\leq C|z| \exp(C|z|)\).
	Plugging the above estimates into \(G\), we get
	\begin{equation*}
		|G(t,y,z)| = \Bigl|\frac{b(t,F(t,y,z))-\partial_{t}F(t,y,z)}{\partial_{y}F(t,y,z)}\Bigr| \leq C(1+|z|)\exp(C|z|).
	\end{equation*}
	One then gets the pathwise estimate
	\begin{equation*}
		\|Y^{n}\|_{\infty} \leq |x_{0}| + T C(1+\|Z^{n}\|_{\infty})\exp(C\|Z^{n}\|_{\infty})\leq C(1+\|Z^{n}\|_{\infty})\exp(C\|Z^{n}\|_{\infty}).
	\end{equation*}
	Again, from the linear growth of \(\sigma\), we have the estimate
	\begin{equation*}
		| F(t,y,z)|\leq C(1+|y|)(1+|z|) \exp(C|z|).
	\end{equation*}
	Therefore, we obtain the pathwise estimate
	\begin{equation*}
		\|X^{n}\|_{\infty}\leq C (1+\|Y^{n}\|_{\infty})(1+\|Z^{n}\|_{\infty}) \exp(C\|Z^{n}\|_{\infty})\leq C(1+\|Z^{n}\|_{\infty})^{2} \exp(2C\|Z^{n}\|_{\infty})\leq C \exp(C\|Z^{n}\|_{\infty}).
	\end{equation*}
	Recall that \(\tilde{Z}^{n}=\E[Z|\mathcal{F}^{n}_{T}]\).
	By applying Jensen's inequality and Fernique's theorem, we deduce the desired uniform moment estimate:
	\begin{equation*}
		\sup_{n}\E[\|X^{n}\|_{\infty}^{\beta}]\leq C \sup_{n}\E[ \exp (C \|Z^{n}\|_{\infty})]\leq C \sup_{n}\E[ \exp (C \|\tilde{Z}^{n}\|_{\infty})]\leq C \E[\exp(C\|Z\|_{\infty})]<\infty.
	\end{equation*}
\end{proof}

\section{Optimality for scalar fractional SDEs}
\label{sec-optim}

We first show two monotonicity properties of the Doss--Sussmann solution map.
\begin{lem}
	\label{lem-mono}
	Assume \(b\in C^{1}_{b}\), \(\sigma\in C^{2}_{b}\), and \(\sigma>0\).
	Let \(z\) be an absolutely continuous path starting at \(0\) and let \(u=\dot{z}\).
	We consider the equation
	\begin{equation*}
		\md x_{t}= b(t,x_{t}) \md t + \sigma(t,x_{t})\md z_{t}.
	\end{equation*}
	Then the solution map \(\Phi(t,x_{0}; u_{\cdot\wedge t})\) is pointwise increasing in \(u\).
	Moreover, if we further assume that
	\begin{equation}
		\label{eqn-mono}
		\partial_{x}b(t,x)\sigma(t,x)- b(t,x) \partial_{x}\sigma(t,x)-\partial_{t}\sigma(t,x)\geq 0,
	\end{equation}
	then the solution map \(\Psi(t,x_{0}; z_{\cdot\wedge t})\) is pointwise increasing in \(z\).
\end{lem}
\begin{rmk}
	We define vector fields \(B=\partial_{t}+ b(t,x) \partial_{x}\) and \(\Sigma= \sigma(t,x) \partial_{x}\).
	The condition \eqref{eqn-mono} is equivalent to \([B,\Sigma]=B \Sigma-\Sigma B\leq 0\).
\end{rmk}
\begin{rmk}
	In the time homogeneous case, the condition \eqref{eqn-mono} is equivalent to \((b/\sigma)\) being increasing.
\end{rmk}
\begin{rmk}
	In the driftless case \(b=0\), the condition \eqref{eqn-mono} is equivalent to \(\sigma(\cdot,x)\) is decreasing for any \(x\).
\end{rmk}
\begin{proof}
	The monotonicity of \(\Phi\) follows immediately from the comparison theorem of ODE.
	Let \(F,G:[0,T]\times \mathbb{R}\times \mathbb{R}\to \mathbb{R}\) be defined as in the previous section.
	We	recall that the solution map \(\Psi\) can be written as
	\begin{equation*}
		\Psi(t,x_{0}; z_{\cdot \wedge t})= F(t, y_{t},z_{t}),
	\end{equation*}
	where  \(y_{\cdot}\) solves
	\begin{equation*}
		\dot{y}_{t}=G(t,y_{t},z_{t}) \text{ with } y_{0}=x_{0}.
	\end{equation*}
	As \(\sigma\geq 0\), we have \(\partial_{z}F\geq 0\).
	Notice that \(H=\partial_{y}F\) satisfies the linear  equation
	\begin{equation*}
		\partial_{z}H(t,y,z)=\partial_{x}\sigma(t,F(t,y,z))H(t,y,z), \text{ with } H(t,y,0)=1.
	\end{equation*}
	Therefore, we deduce \(\partial_{y}F>0\).
	In order to show the monotonicity of \(\Psi\), it suffices to show \(y_{t}\) is pointwise increasing in \(z_{\cdot \wedge t}\).
	A direct computation yields
	\begin{align*}
		\partial_{z}G & = \frac{\partial_{x}b (t,F)\partial_{z}F - \partial_{t}\partial_{z}F }{\partial_{y} F}- \frac{  (b(t,F)-\partial_{t}F)\partial_{yz}F}{(\partial_{y}F)^{2}} \\
		              & = \frac{\partial_{x}b(t,F)\sigma(t,F)- b(t,F) \partial_{x}\sigma(t,F)-\partial_{t}\sigma(t,F)}{\partial_{y}F}\geq 0,
	\end{align*}
	where we plug in \(\partial_{z}F= \sigma(t,F)\) in the second line.
	Therefore, by the comparison theorem of ODE, we deduce \(y_{t}\) is pointwise increasing in \(z_{\cdot \wedge t}\) and we conclude the proof.
\end{proof}

\begin{asmp}
	\label{asmp-submod}
	We assume that \(c\) is submodular in the sense that
	\begin{equation*}
		c(x,y)+c(x',y')\leq c(x,y')+c(x',y)
	\end{equation*}
	for any paths \(x\preceq x'\) and \(y\preceq y'\).
\end{asmp}

\begin{thm}
	\label{thm-opt}
	Suppose that Assumptions \ref{asmp-wp} and \ref{asmp-submod} hold.
	Let \(B\) be a standard Brownian motion on \((\Omega,\mathcal{F},\P)\), and for \(H,\bar{H}\in(0,1)\), let
	\begin{equation*}
		Z_{t}^{H}=\int_{0}^{t}k_{H}(t,s)\md B_{s}, \qquad \bar{Z}_{t}^{\bar{H}}=\int_{0}^{t}k_{\bar{H}}(t,s)\md B_{s}.
	\end{equation*}
	Let \(X,Y\)  be the corresponding Doss--Sussmann solutions with coefficients \((b,\sigma)\) and \((\bar{b},\bar{\sigma})\) respectively.
	We assume either the following conditions hold:
	\begin{enumerate}
		\item the Hurst parameter \(H\in[1/2,1)\);
		\item the Hurst parameter \(H\in(0,1/2)\) and the coefficients \(b,\sigma\) satisfy \eqref{eqn-mono},
	\end{enumerate}
	The analogous assumption also holds for \((\bar{H},\bar{b},\bar{\sigma})\).
	Then the synchronous coupling \((\Id,\Id)_{\#}P\) is the optimal coupling.
\end{thm}

In order to prove the above theorem via our discretization result, we will require the following lemma establishing the optimizers of the discrete time bicausal transport problems:

\begin{lem}
	\label{lem-submodular-backward}
	Let \(\xi=(\xi_{1},\ldots,\xi_{N})\) be a discrete process with natural filtration \(\mathbf{F}\).
	Assume that \((\xi_{i})_{1\leq i\leq N}\) are independent.
	Let \(h:\mathbb{R}^{N}\times\mathbb{R}^{N}\to\mathbb{R}_{\geq0}\) be Borel measurable, submodular, and \(h(x,y)\leq a(x)+b(y)\) for some non-negative \(a,b\in L^{1}(\xi)\).
	Then the synchronous coupling \(\pi_{0}=(\Id,\Id)_{\#}\P\) attains
	\begin{equation*}
		\inf_{\pi\in \Pi_{\bc}(\mathbf{F},\mathbf{F})}\E_{\pi}[h(\xi(\omega),\xi(\bar{\omega}))].
	\end{equation*}
\end{lem}

\begin{proof}[Proof of Theorem \ref{thm-opt}]
	By Theorem \ref{thm-approx}, we have
	\begin{equation*}
		\mathcal{T}_{\bc}(X,Y; c)= \lim_{n\to \infty}\mathcal{T}_{\bc}(\hat{X}^{n},\hat{Y}^{n}; \hat{c}_{n}).
	\end{equation*}
	We claim that the synchronous coupling  is the optimal coupling for the discrete transport problem for any \(n\).
	By Theorem \ref{thm-stab}, it is clear that synchronous coupling is also optimal for the limiting continuous-time bicausal optimal transport problem.

	\emph{Step 1.}
	We first assume \(Z\) is a fractional Brownian motion with \(H\in[1/2,1)\).
	Recall
	\begin{equation*}
		\hat{X}^{n}_{t_{k}}= \Phi(t_{k},x_{0}; \dot{Z}^{n}_{\cdot \wedge t_{k}})
		=\Phi(t_{k},x_{0}; (\Delta Z^{n}_{t_{i}})_{i\leq k-1}),
	\end{equation*}
	where
	\begin{equation*}
		\frac{T}{2^{n}}	\Delta Z^{n}_{t_{i}}=Z^{n}_{t_{i+1}}-Z^{n}_{t_{i}} = \sum_{j=0}^{i-1}\int_{t_{j}}^{t_{j+1}}(k_{H}(t_{i+1},s)-k_{H}(t_{i},s))\md s \Delta B_{t_{j}} + \int_{t_{i}}^{t_{i+1}}k_{H}(t_{i+1},s)\md s \Delta B_{t_{i}}.
	\end{equation*}
	For any \(H\in [1/2,1)\) and \(s\in[0,T]\),  \(k_{H}(\cdot,s)\) is increasing.
	Therefore, together with Lemma \ref{lem-mono}, we deduce \(\hat{X}^{n}_{t_{k}}\) is a \emph{pointwise} increasing functional of \((\Delta B_{t_{i}})_{i\leq k-1}\).

	\emph{Step 2.}
	If we assume \(H\in(0,1/2)\), then
	\begin{equation*}
		\hat{X}^{n}_{t_{k}}=\Psi(t_{k},x_{0};Z^{n}_{\cdot \wedge t_{k} })= \Psi(t_{k},x_{0}; (Z^{n}_{t_{i}})_{i\leq k}).
	\end{equation*}
	We recall that
	\begin{equation*}
		Z^{n}_{t_{i}}=\frac{2^{n}}{T}\sum_{j=0}^{i-1}\int_{t_{j}}^{t_{j+1}}k_{H}(t_{i},s)\md s\Delta B_{t_{j}}.
	\end{equation*}
	Therefore, by Lemma \ref{lem-mono} and \(k_{H}(t,s)\geq 0\), we again deduce \(\hat{X}^{n}_{t_{k}}\) is a \emph{pointwise} increasing functional of  \((\Delta B_{t_{i}})_{0\leq i\leq k-1}\).

	\emph{Step 3.}
	Following the same argument, \(\hat{Y}^{n}_{t_{k}}\) is also a \emph{pointwise} increasing functional of \((\Delta B_{t_{i}})_{0\leq i\leq k-1}\) for any \(H\in(0,1)\).
	Therefore,   \(\hat{X}^{n}= \Gamma_{n}^{X}(\Delta B)\) and \(\hat{Y}^{n}=\Gamma_{n}^{Y}(\Delta B)\) with \(\Gamma_{n}^{X},\Gamma_{n}^{Y}: \mathbb{R}^{|D_{n}|}\to \mathbb{R}^{|D_{n}|}\) triangular and pointwise increasing.
	Now we consider the discrete-time bicausal optimal transport problem \(\mathcal{T}_{\bc}(\hat{X}^{n},\hat{Y}^{n}; \hat{c}_{n})\).
	Recall we set \(\hat{c}_{n}=c(l_{n},l_{n})\).
	We write \(h:\mathbb{R}^{|D_{n}|}\times \mathbb{R}^{|D_{n}|}\to \mathbb{R}\) as
	\begin{align*}
		\hat{c}_{n}(\hat{X}^{n}(\omega),\hat{Y}^{n}(\bar{\omega})) & = c(l_{n}\circ \Gamma_{n}^{X}(\Delta B (\omega)),l_{n}\circ \Gamma_{n}^{Y}(\Delta B (\bar{\omega}))) \\
		                                                           & :=h( \Delta B(\omega), \Delta B(\bar{\omega})).
	\end{align*}
	Since \(c\) is submodular and \(l_{n}\circ \Gamma_{n}^{X}(x)\preceq l_{n}\circ \Gamma_{n}^{X}(y)\) for any  \(x\preceq y\), we deduce that \(h\) is submodular.
	We apply Lemma \ref{lem-submodular-backward} with \(N=|D_{n}|\), \(\xi_{i}= \Delta B_{t_{i}}\) and the submodular cost \(h\) given above.
	This shows that the synchronous coupling is optimal for the discrete-time problem with any dyadic discretization.
	Passing \(n\to \infty\), following the proof of Theorem \ref{thm-stab}, we derive the optimality of the synchronous coupling for the original continuous-time problem.
\end{proof}

\begin{proof}[Proof of Lemma \ref{lem-submodular-backward}]
	We write \(\xi_{\#}P=\lambda\) and \((\xi_{k})_{\#}P=\lambda_{k}\) for \(1\leq k\leq N\).
	For \(0\leq k\leq N\), define
	\begin{align*}
		h_{k}(x_{1:k},y_{1:k}):=\int_{\mathbb{R}^{N-k}} & h((x_{1:k},z_{k+1:N}),(y_{1:k},z_{k+1:N})) \lambda_{k+1}(\md z_{k+1})\otimes \cdots \otimes \lambda_{N}(\md z_{N}),
	\end{align*}
	with the conventions \(h_{N}=h\) and
	\(
	h_{0}=\int_{\mathbb{R}^{N}}h(z,z)\lambda(\md z).
	\)
	We note every \(h_{k}\) is submodular since	\((x_{1:k},y_{1:k})\mapsto h(x_{1:k},z_{k+1:N}, y_{1:k},z_{k+1:N})\) is submodular for  any  fixed \(z_{k+1:N}\).
	Fix \(\pi\in\Pi_{\bc}(\mathbf{F},\mathbf{F})\) and consider the discrete time stochastic processes $ \Omega \times \Omega \ni (\omega, \overline{\omega}) \mapsto \xi(\omega)$ and $(\omega, \overline{\omega}) \mapsto \overline{\xi}(\omega, \overline{\omega}):= \xi(\overline{\omega})$.
	By \cite[Lemma A.1]{bartl24Wasserstein}, conditionally on \((\xi_{1:N-1},\overline{\xi}_{1:N-1})\), the pair \((\xi_{N},\overline{\xi}_{N})\) has a joint law with both marginals given by \(\lambda_{N}\).
	Since \(h\) is submodular, we note
	\begin{equation*}
		(u,v)\mapsto g(u,v)=h_{N}((\xi_{1:N-1}^{1},u),(\xi_{1:N-1}^{2},v))
	\end{equation*}
	is submodular.
	In particular, the diagonal set \(\{(u,u): u\in \mathbb{R}\}\) is a \(g\)-monotone set.
	Therefore, we have \(\pi\)-a.s.
	\begin{align*}
		\E_{\pi}[h_{N}(\xi_{1:N},\overline{\xi}_{1:N})| (\xi_{1:N-1},\overline{\xi}_{1:N-1})] & \geq\int_{\mathbb{R}}h_{N}((\xi_{1:N-1},u),(\overline{\xi}_{1:N-1},u))\lambda_{N}(\md u) \\
		                                                                                      & =h_{N-1}(\xi_{1:N-1},\overline{\xi}_{1:N-1}).
	\end{align*}
	By the tower law of conditional expectation, we derive
	\begin{equation*}
		\E_{\pi}[h_{N}(\xi_{1:N},\overline{\xi}_{1:N})]\geq \E_{\pi}[h_{N-1}(\xi_{1:N-1},\overline{\xi}_{1:N-1})].
	\end{equation*}
	Repeating the same argument for \(k=N-1,\ldots,1\) yields
	\begin{equation*}
		\E_{\pi}[h(\xi,\overline{\xi})]\geq h_{0}=\E_{\pi_{0}}[h(\xi,\xi)].
	\end{equation*}
	We conclude the proof by noticing that the synchronous coupling \(\pi_{0}\) is bicausal.
\end{proof}

A final application of our stability results allows us to recover synchronous optimality for It\^o SDEs under mild assumptions on the coefficients as in \citet{backhoff-veraguas22Adapted}.

\begin{cor}
	\label{cor:sdesyncoptimal}
	Let $B$ be a standard Brownian motion on $(\Omega, \mathcal{F}, P)$ and $X,Y$ be the pathwise unique strong solution to the \emph{It\^o} stochastic differential equations \[
		dX_t = b(t,X_t) dt + \sigma(t,X_t) dB_t, \quad X_0 = x_0 \text{ and }     dY_t = \overline{b}(t,Y_t) dt + \overline{\sigma}(t,Y_t) dB_t, \quad Y_0 = y_0
	\] respectively. We assume that the cost $c$ satisfies Assumption \ref{asmp-wp} $(iii)$ and \ref{asmp-submod} and either of the following conditions hold: \begin{enumerate}[(i)]
		\item $b,\sigma$ are continuous with linear growth with $\sigma \geq 0$;
		\item $b,\sigma$ are bounded, measurable and $\sigma > c$ for some $c > 0.$
	\end{enumerate} The analogous assumption also holds for $(\overline{b}, \overline{\sigma})$. Then, the synchronous coupling $(\Id,\Id)_{\#}P$ is the optimal coupling.
\end{cor}

\begin{rmk}
	\label{rmk-rightctsfiltrationweakuniqueness}
	Note that we implicitly assume the existence of a pathwise unique strong solution and by the Yamada-Watanabe theorem, this gives weak uniqueness, i.e. uniqueness for the martingale problem. By \cite[Theorem $4.29$, Proposition $4.32$]{jacod2003limit}, the (completed) filtration generated by $X$ and $Y$ are right continuous.
\end{rmk}

\begin{proof}
	Up to standard mollification and truncation techniques, there exist sequences $b_n, \sigma_n \in C^\infty_b$ such that if (i) holds for $(b,\sigma)$, then \[
		b_n,\sigma_n \text{ converges uniformly on compacts to } b, \sigma \text{ and } \sup_{n,t,x} \frac{|\sigma_n(t,x)| + |b_n(t,x)|}{1 + |x|}<\infty \text{ and } \sigma_n > 0
	\] or if (ii) holds, then \[
		b_n,\sigma_n \text{ converges in $L^3_{loc}$ to } b, \sigma \text{ and } \sup_{n,t,x} |\sigma_n(t,x)|+|b_n(t,x)| < \infty \text{ and } \inf_n \sigma_n > \delta
	\] for some $\delta > 0.$ Since $b_n, \sigma_n \in C^\infty_b$, then we consider their pathwise unique strong solutions $X^n,Y^n$ to the It\^{o} stochastic differential equation driven by Brownian motion $B$ \[
		dX^n = b_n(t,X^n_t) dt + \sigma_n(t,X^n_t) dB_t, \quad X^n_0 = x_0
	\] with Stratonovich form \[
		dX^n = (b_n(t,X^n_t) - \frac{1}{2} \sigma_n(t,X^n_t) \partial_x \sigma_n(t,X^n_t)) dt + \sigma_n(t,X^n_t) \circ dB_t, \quad X^n_0 = x_0
	\] and a similar equation holding for $Y^n$. Then, Theorem \ref{thm-opt} yields that the synchronous coupling is optimal for $\mathcal{T}_{\bc}(X^n,Y^n; c)$ for every $n \geq 1$, i.e. \[
		\mathcal{T}_{\bc}(X^n,Y^n; c) = E_P [c(X^n,Y^n)].
	\] We claim that as $n \rightarrow \infty$, the LHS and RHS converges to $\mathcal{T}_{\bc}(X,Y;c)$ and $E_P[ c(X,Y)]$ respectively, which yields the corollary. We note that using Theorem $A$ and $B$ in \cite{kaneko1988note} under condition (i) and (ii) respectively, we obtain \begin{equation}
		\label{eqn-sdecvg}
		\lim_{n\to\infty} E_P[ \Vert X^n - X \Vert_\infty^2] + E_P [\Vert Y^n - Y \Vert_\infty^2] =0.
	\end{equation} Furthermore, since $b_n,\sigma_n$ are linear growth with linear growth constant independent of $n$, then $X^n,Y^n, n \geq 0$ satisfies Assumption \ref{asmp-conv} $(iii)$ for any $\beta > 0$. In particular, by Assumption \ref{asmp-wp} $(iii)$ on the cost $c$, uniform integrability and Vitali's convergence theorem, we obtain \[
		\lim_{n\to\infty} E_P [c(X^n,Y^n)] = E_P [c(X,Y)].
	\] To show the convergence of the LHS, we verify the various conditions of Theorem \ref{thm-stab}. Since $X^n,Y^n,X$ and $Y$ are solutions to SDEs with pathwise unique strong solutions driven by the \emph{same} Brownian motion $B$, \cite[Corollary $3.3$]{cont2024causal} yields that $(\Id,\Id)_{\#}P \in \Pi_{\bc}(\mathbf{F},\mathbf{F}^n)$ and hence, \[
		\lim_{n\to\infty} \mathcal{A} \mathcal{W}(X^n, X) + \mathcal{A} \mathcal{W}(Y^n, Y)= 0
	\] and condition $(i)$ is satisfied with $\bt = (\Id,\Id,\dots)$. Condition $(ii)$ is satisfied by the existence of pathwise unique strong solutions (see Remark \ref{rmk-rightctsfiltrationweakuniqueness}). Condition $(iii)$ is trivially satisfied. Condition $(iv)$ follows from the fact that $c$ satisfies Assumption \ref{asmp-wp} $(iii)$ and $X^n,Y^n, n \geq 0$ satisfies Assumption \ref{asmp-conv} $(iii)$ for any $\beta > 0$. Hence, Theorem \ref{thm-stab} yields \[
		\lim_{n\to\infty} \mathcal{T}_{\bc}(X^n,Y^n; c) = \mathcal{T}_{\bc}(X,Y;c).
	\]
\end{proof}

\bibliography{mybib}

\end{document}